\documentclass[hidelinks,onefignum,onetabnum]{siamart251216}

\usepackage{lipsum}
\usepackage{amsfonts}
\usepackage{graphicx}
\usepackage{epstopdf}
\usepackage{algorithmic}
\usepackage{amsopn}
\usepackage{booktabs}
\usepackage{amsmath, amssymb,mathrsfs, caption, subcaption}
\usepackage{bm}
\usepackage{tikz}
\usepackage{multirow}
\usepackage{siunitx}
\usepackage{comment}

\newcommand{\dx}{\;{\rm d}\bm x}
\newcommand{\ds}{\;{\rm d} s}
\ifpdf
  \DeclareGraphicsExtensions{.eps,.pdf,.png,.jpg}
\else
  \DeclareGraphicsExtensions{.eps}
\fi

\newsiamremark{hypothesis}{Hypothesis}
\crefname{hypothesis}{Hypothesis}{Hypotheses}
\newsiamthm{claim}{Claim}
\newsiamremark{fact}{Fact}
\crefname{fact}{Fact}{Facts}

\headers{An Incidence-Based Edge-Averaged Method}{J. H. Adler, X. Hu, S. Lee, and E. O'Keefe}

\title{An Incidence-Based Edge-Averaged Method for Convection-Diffusion Problems on Polytopal Meshes\thanks{Submitted to the editors DATE.}}

\author{James H. Adler\thanks{Department of Mathematics, Tufts University, Medford, MA 
  (\email{james.adler@tufts.edu}, \email{xiaozhe.hu@tufts.edu}, \email{seulip.lee@tufts.edu},
  \email{eoghan.okeefe@tufts.edu}).}
  \and Xiaozhe Hu\footnotemark[2]
  \and Seulip Lee\footnotemark[2]
  \and Eoghan O'Keefe\footnotemark[2]
}

\usepackage{amsopn}

\usepackage{nicematrix}

\ifpdf
\hypersetup{
  pdftitle={An Example Article},
  pdfauthor={D. Doe, P. T. Frank, and J. E. Smith}
}
\fi

\begin{document}

\maketitle

\begin{abstract}
We develop an incidence-based edge-averaged method for scalar convection-diffusion problems on convex polytopal meshes in arbitrary dimensions. The proposed method extends the classical edge-averaged finite-element scheme through a graph-theoretic formulation based on vertex-edge incidence matrices. We first establish this formulation on simplices in arbitrary dimensions and then extend it to convex polytopes by introducing complementary edges and generalized barycentric coordinates. A key ingredient is a matrix identity that relates the generalized Whitney edge mass matrix to a diagonal representation of the Poisson stiffness matrix, enabling an exponentially-fitted discretization expressed entirely in terms of nodal values and edge weights. We prove that the resulting stiffness matrix is a nonsingular $M$-matrix under suitable positivity and connectivity conditions, ensuring monotonicity and the discrete maximum principle. Numerical experiments on two- and three-dimensional simplicial and nonsimplicial meshes demonstrate first-order convergence in the $H^1$ semi-norm and stable, nonoscillatory solutions in convection-dominated regimes.
\end{abstract}

\begin{keywords}
convection-diffusion equations; edge-averaged method; incidence matrix; polytopal meshes; arbitrary dimensions; monotonicity.
\end{keywords}

\begin{MSCcodes}
65N30, 65N12, 05C50
\end{MSCcodes}

\section{Introduction}
Convection-diffusion equations arise in many applications involving the combined effects of transport and diffusion. A major numerical difficulty occurs in the convection-dominated regime, where solutions may develop sharp boundary or interior layers. Unless the mesh is sufficiently refined and suitably aligned with these layers, standard Galerkin finite-element discretizations often produce nonphysical oscillations that substantially deteriorate the quality of the numerical solution \cite{JOHN20072197,1050003824962563456,Stynes_2005}. This difficulty has motivated the development of stabilized discretizations for convection-dominated problems \cite{AUGUSTIN20113395,survey}.

Methods to address such issues have been developed mainly in two directions. The first introduces additional stabilization terms into the weak formulation to obtain stability and convergence estimates that remain robust as the diffusion coefficient approaches zero (see \cite{BROOKS1982199,BURMAN20041437,lps} for examples). Although these methods can reduce spurious oscillations, they generally do not guarantee a discrete maximum principle (DMP), which ensures the discrete solution is free of nonphysical over- and undershoots. Hence, the second direction seeks to enforce the DMP. Algebraic flux correction methods \cite{barrenechea2018unified} do so directly, but typically lead to nonlinear algebraic systems. Linear monotone discretizations offer another approach, as monotonicity is a sufficient condition for the discrete maximum principle. At the algebraic level, a stronger but convenient sufficient condition for monotonicity is that the stiffness matrix be a nonsingular \(M\)-matrix, which can be checked directly from the assembled matrix.

Another monotone approach, the edge-averaged finite element (EAFE) method \cite{eafe}, has been developed for convection-diffusion problems on simplicial meshes. Motivated by the Scharfetter-Gummel scheme \cite{original_Scharfetter-Gummel_paper}, EAFE constructs an edgewise approximation of the exponentially fitted flux. Related edge-averaging and exponential-fitting ideas have subsequently been used to develop a variety of numerical schemes \cite{bank1998finite,brezzi1989two,dorfler1999uniform,Scharfetter-Gummel_extensions}. The monotonicity of EAFE is guaranteed when the stiffness matrix for the Poisson equation in the same conforming finite-element space is an \(M\)-matrix, a condition that is satisfied, for example, on Delaunay triangulations in two dimensions. An important feature of EAFE is that the final formulation can be expressed in terms of nodal differences and diffusion stiffness weights without explicitly computing the associated edge basis functions.

Several recent works have extended different aspects of EAFE. The simplex-averaged finite element (SAFE) method \cite{safe} generalizes the approach from scalar \(H(\mathrm{grad})\) problems to vector convection-diffusion problems in \(H(\mathrm{curl})\) and \(H(\mathrm{div})\), using graph coboundary operators on simplicial meshes. A related stable mimetic finite-difference (MFD) method \cite{casey2} employs primal-dual meshes, incidence-relation matrices, discrete differential operators, and shifted flux operators to obtain scalar and vector schemes in a finite-difference setting. The edge-averaged virtual element (EAVE) method \cite{eave}, on the other hand, extends EAFE to scalar convection-diffusion problems on polygonal meshes in two dimensions. It includes both a monotone formulation on suitable Voronoi meshes and a general formulation on arbitrary polygonal meshes. However, these developments do not yet provide a graph-based, edge-averaged method for general polytopal meshes in three or more dimensions. Both SAFE and MFD have primarily been developed on simplicial meshes in two and three dimensions, whereas the existing EAVE construction is limited to two-dimensional polygonal meshes. These limitations are particularly relevant for convection-dominated problems, where adaptive or anisotropic refinement and layer-aligned mesh generation may benefit from greater mesh flexibility and the use of non-simplicial elements. This need becomes more pronounced in three-dimensional problems, where constructing high-quality simplicial meshes can be challenging, and in space-time discretizations \cite{bank2017arbitrary}, where treating time as an additional coordinate may require four-dimensional meshes for time-dependent three-dimensional problems \cite{adler2025unified}.

Therefore, in this work, we propose an incidence-based edge-averaged method (IBEAM) for scalar convection-diffusion problems on polytopal meshes in arbitrary dimensions. We first develop an incidence-based formulation of EAFE on simplices in arbitrary dimensions and then extend this framework to general convex polytopes.
The key observation underlying our construction is that the one-skeleton of a simplex is a complete graph, whereas that of a general polytope is not. We therefore augment the local edge graph of each polytopal element to the complete graph on its vertex set. The additional connections, which correspond geometrically to element diagonals, are referred to as complementary edges and are used as auxiliary local objects in the construction. On the resulting complete graph, we define an extended incidence matrix that includes both geometric and complementary edges. We also consider a nodal finite-element space generated by generalized barycentric coordinates, primarily Wachspress coordinates \cite{wachspress3d,wachspress2011barycentric} in this work, together with an associated edge space constructed from Whitney one-forms \cite{lohi2021whitney} on the extended edge set. The extended incidence matrix represents the gradients of the nodal basis functions in terms of these edge functions, thereby connecting the nodal stiffness matrix to the mass matrix of the extended edge space. We then replace the generally dense edge mass matrix by a diagonal representation that retains the required relation with the nodal stiffness matrix. Combining this representation with an exponentially-fitted flux operator yields the IBEAM discretization without requiring the explicit computation of the extended edge basis functions. Since its derivation is based on incidence relations between compatible local vertex and edge spaces, the algebraic construction is not intrinsically restricted to a particular spatial dimension and can be extended to higher dimensions. As in EAFE, when the stiffness matrix of the associated nodal discretization for the Poisson equation is an \(M\)-matrix, so is the resulting IBEAM stiffness matrix, providing a sufficient condition for monotonicity. Numerical experiments on simplicial and nonsimplicial meshes demonstrate first-order convergence in the $H^1$ semi-norm and stable, nonoscillatory behavior in convection-dominated regimes.

The remainder of this paper is organized as follows. \Cref{section:prelim} introduces the convection-diffusion model problem and reviews the classical EAFE scheme.
\Cref{sec:graph} develops an incidence-based formulation on simplices in arbitrary dimensions and extends the framework to general convex polytopes using complementary edges and generalized barycentric coordinates.
\Cref{section:scheme} presents the IBEAM discretization and establishes sufficient conditions for its monotonicity.
\Cref{section:numerical} investigates the convergence behavior and numerical performance of IBEAM on various simplicial and polytopal meshes, including convection-dominated problems with boundary layers. Finally, \Cref{section:conclusion} summarizes our contributions and discusses future research directions.

\section{Preliminaries}\label{section:prelim}

This section introduces the convection-diffusion model problem and reviews the classical EAFE scheme \cite{eafe}, providing the necessary background for the graph-theoretic framework developed in the next section.

\subsection{Model Problem}

Let $\Omega\subset\mathbb{R}^d$ be a bounded Lipschitz domain. We denote by $H^1(\Omega)$ and $H_0^1(\Omega)$ the standard Sobolev spaces, and by $L^2(\Omega)$ the space of square-integrable functions.
The norm, semi-norm, and inner product associated with a Hilbert space $X$ are denoted by $\|\cdot\|_X$, $|\cdot|_X$, and $\langle\cdot,\cdot\rangle_X$, respectively.

We consider the following steady-state convection-diffusion problem:
\begin{subequations}\label{sys: conv_diff_eqn}
\begin{alignat}{2}
-\nabla \cdot \left(\alpha \nabla u + \bm \beta u\right) &= f \quad &&\text{in } \Omega, \\
u &= 0 \quad &&\text{on }\partial \Omega,
\end{alignat}
\end{subequations}
where $f \in L^2(\Omega)$, $\alpha \in C^0(\Bar{\Omega})$ is the diffusion coefficient satisfying
\[
0 < \alpha_{\min}\leq \alpha(\bm x)\leq \alpha_{\max}
\qquad \forall \bm x\in\bar{\Omega},
\]
and $\bm \beta \in [C^0(\Bar{\Omega})]^d$ is the convection field.
We are particularly interested in convection-dominated regimes, where the magnitude of the convective transport is large relative to diffusion. Defining the flux by
\[
    J(u)=\alpha\nabla u+\bm\beta u,
\]
the weak formulation is to find $u\in H_0^1(\Omega)$ such that
\[
    B(u,v)=F(v) \qquad \forall v\in H_0^1(\Omega),
\]
where 
\begin{equation}\label{bilinear_form}
    B(u,v):=\int_\Omega J(u)\cdot\nabla v\dx,
    \qquad
    F(v):=\int_\Omega fv\dx.
\end{equation}

\subsection{Finite-Element Spaces}

We introduce the finite-element spaces and notation underlying the classical EAFE scheme. Although EAFE is formulated in two and three dimensions, we restrict our discussion to triangular meshes for simplicity. In particular, we describe the linear Lagrange and lowest-order N\'ed\'elec finite-element spaces and the relationship between their basis functions.

Let $\mathcal{T}_h$ be a conforming triangulation of $\Omega$, with vertex set $\mathcal{V}_h$ and edge set $\mathcal{E}_h$. 
For each triangle $T\in\mathcal{T}_h$, let $\{\bm x_i\}_{i=1}^3$ denote its vertices and $\{e_l\}_{l=1}^3$ its edges.
We orient each edge $e_l = (i,j)$ from the lower-indexed vertex, $\bm x_i$ to the higher-indexed vertex, $\bm{x}_j$ and define the corresponding oriented edge vector by

\begin{equation}\label{eq:tangent_vector}
    \bm\tau_l=
    \bm x_j-\bm x_i,\qquad i<j.
\end{equation}

We define the conforming finite-element spaces
\begin{align*}
    V_h
    &=
    \{v_h\in H_0^1(\Omega)\,:\,
    v_h|_T\in \mathcal{P}_1(T),
    \ \forall T\in\mathcal{T}_h\},\\
    Q_h
    &=
    \{\bm w_h\in H(\mathrm{curl};\Omega)\,:\,
    \bm w_h|_T\in\mathcal N_0(T),
    \ \forall T\in\mathcal{T}_h\},
\end{align*}
where $\mathcal{P}_1(T)$ denotes the space of linear polynomials on $T$, and $\mathcal N_0(T)$ denotes the local lowest-order N\'ed\'elec space, given by
\[
    \mathcal{N}_0(T) = [\mathcal{P}_0(T)]^2 + \bm x^\perp \mathcal{P}_0(T),
\]
where $\bm x^\perp$ is a $90^\circ$ rotation of $\bm x$ satisfying $\bm x \cdot \bm x^\perp = 0$.
The degrees of freedom for any $\bm w \in \mathcal{N}_0(T)$ are defined by
\[
    \text{dof}_{l}(\bm w) = \frac{1}{|e_l|} \int_{e_l} \bm w\cdot \bm \tau_l\ds = \left.(\bm w\cdot \bm \tau_l)\right|_{e_l},
\]
where $e_l = (i, j)$.

Let $\{\lambda_i\}_{i=1}^3$ denote the barycentric coordinates on $T$, with the local vertices $\{\bm x_i\}_{i=1}^3$ of $T$ numbered in increasing order of their global indices, and let $\{\bm \phi_l\}_{l=1}^3$ denote the local canonical N\'ed\'elec basis functions satisfying $\mathrm{dof}_l(\bm \phi_m) = \delta_{lm}$.
For each edge $e_l = (i,j)$, these basis functions are given by
\[
    \bm \phi_l = \lambda_i \nabla \lambda_j - \lambda_j \nabla \lambda_i.
\]

\subsection{Edge-Averaged Finite-Element Scheme}\label{subsec: EAFE}

We now review the classical edge-averaged finite-element (EAFE) scheme \cite{eafe}, focusing on the key steps in its derivation that are relevant to the subsequent generalization.

For each triangle $T\in \mathcal{T}_h$, we define the local bilinear form
\[
    B_T(u_h,v_h) = \int_T J(u_h) \cdot \nabla v_h\dx,
\]
where \(u_h,v_h \in V_h\).
We approximate the flux \(J(u_h)\) by an exponentially-fitted flux $J_T(u_h)\in \mathcal{N}_0(T)$, yielding
\begin{equation*}
B_T^\mathcal{N}(u_h,v_h) = \int_T J_T(u_h)\cdot \nabla v_h\dx,
\end{equation*}
where
\begin{equation*}
J_T(u_h)
= \sum_{l=1}^3 \mathrm{dof}_l(J_T(u_h))\bm \phi_l
= \sum_{l=1}^3
\left.(J_T(u_h)\cdot \bm \tau_l)\right|_{e_l}\bm \phi_l,
\end{equation*}
and \(J_T(u_h)\in\mathcal{N}_0(T)\).
Since \(\left.\nabla v_h\right|_T\in \mathcal{N}_0(T)\), we express the local bilinear form as
\begin{align*}
B_T^\mathcal{N}(u_h,v_h)
&= \int_T J_T(u_h)\cdot \nabla v_h\dx \\
&= \sum_{l=1}^3\sum_{m=1}^3
\left.(J_T(u_h)\cdot \bm \tau_l)\right|_{e_l}
\left(\nabla v_h\cdot \bm \tau_m)\right|_{e_m}
\int_T \bm \phi_l \cdot \bm \phi_m\dx.
\end{align*}
To express these edge quantities in terms of nodal values, we introduce the edge difference operator
\begin{align*}
\delta_l u_h := u_h(\bm x_j)-u_h(\bm x_i)
= u_j-u_i,
\end{align*}
where \(e_l=(i,j)\) is edge $l$ connecting vertex $i$ to vertex $j$. 
Since $v_h$ is linear on $T$, the gradient term satisfies 
\[\left.(\nabla v_h\cdot\bm\tau_m)\right|_{e_m} = \delta_m v_h.
\]
Next, we consider the flux term. For each edge \(e_l\), we introduce a locally defined potential function \(\gamma\) satisfying 
\[\nabla\gamma\cdot\bm\tau_l = \alpha^{-1}\bm\beta\cdot\bm\tau_l\qquad \text{on }e_l.
\]
From the definition of the flux, we have
\begin{align*}
\alpha^{-1}e^\gamma(J(u_h)\cdot\bm\tau_l)
=
\nabla(e^\gamma u_h)\cdot
\bm\tau_l.
\end{align*}
Applying $\mathrm{dof}_l$ to both sides through the weighted flux interpolation and using the fact that the tangential component of \(J_T(u_h)\in\mathcal N_0(T)\) is constant along \(e_l\), we obtain

\[
\left.(J_T(u_h)\cdot\bm\tau_l)\right|_{e_l}
\frac{1}{|e_l|}
\int_{e_l}\alpha^{-1}e^\gamma\dx  =
\frac{1}{|e_l|}\int_{e_l}
\nabla(e^\gamma u_h)\cdot
\bm\tau_l\dx
=
\delta_l(e^\gamma u_h).
\]
Therefore,
\begin{equation}
\left.(J_T(u_h)\cdot\bm\tau_l)\right|_{e_l}
=
\kappa_l\delta_l(e^\gamma u_h),\label{eqn: exponential_fitting}
\end{equation}
where
\begin{align*}
\kappa_l
=
\left(
\frac{1}{|e_l|}
\int_{e_l}\alpha^{-1}e^\gamma\dx
\right)^{-1}.
\end{align*}
Substituting this expression and the corresponding edge difference for \(v_h\) into the local bilinear form yields
\begin{equation}\label{eqn: local_B^N}
B_T^\mathcal{N}(u_h,v_h)
=
\sum_{l=1}^3\sum_{m=1}^3
\kappa_l\delta_l(e^\gamma u_h)
\delta_m v_h
\int_T\bm\phi_l\cdot\bm\phi_m\dx.
\end{equation}
To derive the EAFE scheme, we employ a mass-lumping technique for the N\'ed\'elec mass matrix that preserves the following identity for the Poisson bilinear form in the linear Lagrange finite-element space \cite{haugazeau1993condensation}:
\begin{equation}\label{eq:simplification_poisson}
\int_T\nabla u_h\cdot\nabla v_h\dx
=
\sum_{l=1}^3\sum_{m=1}^3
\delta_lu_h\delta_mv_h
\int_T\bm\phi_l\cdot\bm\phi_m\dx
=
\sum_{l=1}^3\omega_l^T\delta_lu_h\delta_lv_h,
\end{equation}
where
\begin{equation}\label{eqn: weights}
\omega_l^T
=
-\int_T\nabla\lambda_i\cdot\nabla\lambda_j\dx,
\end{equation}
for the edge \(e_l=(i,j)\).
Motivated by this identity, we approximate the double-sum formulation in \eqref{eqn: local_B^N} by replacing the N\'ed\'elec mass matrix with the edge weights $\omega_l^T$, leading to the local EAFE bilinear form
\begin{equation}\label{eafe_bilinear_form}
B_T^h(u_h,v_h)
=
\sum_{l=1}^3
\omega_l^T\kappa_l
\delta_l(e^\gamma u_h)\delta_lv_h.
\end{equation}
We note that when $J_T(u_h)$ is constant on $T$, it can be expressed as $J_T(u_h) = \nabla\psi$ for some linear function $\psi$. In this case, the Poisson identity implies
\begin{equation}
    B_T^\mathcal{N}(u_h,v_h) = B_T^h(u_h,v_h).\label{eqn: equiv_B^N_B^h}
\end{equation}
For a nonconstant $J_T(u_h)$, \(B_T^h\) instead provides an approximation to \(B_T^{\mathcal N}\).
The EAFE scheme has been shown to achieve first-order convergence with respect to the mesh size $h$, under suitable regularity assumptions and for sufficiently small $h$ \cite{eafe}. The convergence analysis relies on a discrete inf-sup condition and a consistency estimate, the latter of which exploits the identity \eqref{eqn: equiv_B^N_B^h} for constant interpolated fluxes.

The EAFE scheme also guarantees monotonicity under the Delaunay mesh condition, preventing spurious oscillations in convection-dominated regimes.
Moreover, its discrete bilinear form involves only edge differences and the corresponding weights, allowing for a simple implementation using the Bernoulli function \cite{eave}.
Although originally formulated on simplicial meshes in two and three dimensions, the EAFE scheme can be naturally interpreted from a graph-theoretic perspective. This perspective not only provides a better understanding of the underlying algebraic structure of EAFE but also enables its extension to simplices in arbitrary dimensions. Furthermore, it provides a systematic approach to deriving generalized EAFE schemes for convection-diffusion problems on non-simplicial polytopal meshes. 

\section{Graph-Theoretic Framework}\label{sec:graph}
This section introduces a graph-theoretic representation of the EAFE scheme and extends this representation to polytopal elements through their embedding into a complete graph. This representation extends the edge-based structure underlying EAFE by introducing complementary edges, which allow the diagonal representation of the Poisson stiffness matrix \eqref{eq:simplification_poisson} to be generalized from simplices to arbitrary convex polytopes.

\subsection{Graph-Based EAFE in Arbitrary Dimensions}\label{subsec: graph-based_EAFE}
We begin by considering a $d$-dimensional simplex $T\in\mathcal{T}_h$ with vertices $\{\bm x_i\}_{i=1}^{n_v}\subset \mathbb{R}^d$ and edges $\{e_l\}_{l=1}^{n_e}$, where $n_v = d+1$ and $n_e = d(d+1)/2$. Each edge is oriented from the lower-indexed vertex to the higher-indexed vertex, consistent with the definition of $\bm \tau_l$ in \eqref{eq:tangent_vector}.

Let $\{\lambda_i\}_{i=1}^{n_v}$ denote the barycentric coordinates on $T$, satisfying
\[
\lambda_i(\bm x_j) = \delta_{ij},\qquad \sum_{i=1}^{n_v}\lambda_i=1.
\]
For each oriented edge \(e_l=(i,j)\), which connects vertex $i$ to vertex $j$, we define the corresponding canonical Whitney edge basis function by
\[
\bm\phi_l = \lambda_i\nabla\lambda_j - \lambda_j\nabla\lambda_i,\qquad i<j.
\]
We define the vertex-edge incidence matrix $G \in \mathbb{R}^{n_e \times n_v}$ by
\[
    G_{lk} = \begin{cases}
        1, & \mathrm{if}\; k=j, \\
        -1, & \mathrm{if}\; k=i, \\
        0, & \mathrm{otherwise,}
    \end{cases}
\]
for each oriented edge $e_l = (i,j)$. The incidence matrix represents the discrete gradient operator, mapping nodal values to their differences along oriented edges. In particular,
\[
\delta_lu_h = (G\bm u_h)_{l},
\]
where $\bm u_h$ denotes the vector of nodal values of $u_h$.
Moreover, the partition-of-unity property of the barycentric coordinates implies
\[
\sum_{j=1}^{n_v} \nabla \lambda_j = \mathbf{0}.
\]
Together with the definition of the Whitney edge functions, this yields (see, e.g., \cite{witney})
\begin{equation}\label{eq:edge_relations}
\nabla \lambda_i
= \sum_{l=1}^{n_e} G_{li}\bm \phi_l.
\end{equation}
This identity expresses the gradients of the nodal basis functions in terms of the Whitney edge functions through the incidence matrix.
Using this relation, the Poisson bilinear form can be expressed as (see \cite{maxwell}) 
\begin{equation}
    \int_T \nabla u_h \cdot \nabla v_h \dx
= \sum_{l=1}^{n_e}\sum_{m=1}^{n_e} \delta_l u_h \delta_m v_h\int_T \bm \phi_l \cdot \bm \phi_m\dx = \bm v_h^\top G^\top M G\bm u_h,\label{eqn: GMG_simplex}
\end{equation}
where \(M\in\mathbb R^{n_e\times n_e}\) is the edge mass matrix defined by
\[
M_{lm} = \int_T \bm \phi_l \cdot \bm \phi_m\dx.
\]
On the other hand, the Poisson bilinear form admits the following representation in terms of edge differences \cite{eafe}:
\begin{equation}
    \int_T \nabla u_h \cdot \nabla v_h \dx = \sum_{l=1}^{n_e} \omega_l^T \delta_l u_h\delta_lv_h = \bm v_h^\top G^\top W G\bm u_h,\label{eqn: GWG_simplex}
\end{equation}
where
\[
W = \mathrm{diag}\left(\omega_1^T,\dots, \omega_{n_e}^T\right),
\]
with the edge weights \(\omega_l^T\) defined as in \eqref{eqn: weights}, now on a \(d\)-dimensional simplex. Since both representations \eqref{eqn: GMG_simplex} and \eqref{eqn: GWG_simplex} hold for arbitrary nodal vectors \(\bm u_h\) and \(\bm v_h\), we obtain the matrix identity,

\begin{equation}\label{g-identity}
G^\top MG = G^\top W G,
\end{equation}
which we use later. A corresponding identity for diagonal mass lumping on tetrahedral meshes can also be found in \cite{bossavit1999yee}.

We now express the EAFE bilinear form \eqref{eafe_bilinear_form} in terms of the graph incidence matrix. Recall that the edge differences of $v_h$ satisfy $\delta_lv_h=(G\bm v_h)_l$. Moreover, following the same argument used to derive \eqref{eqn: exponential_fitting},
\[
    \left.(J_T(u_h)\cdot \bm \tau_l)\right|_{e_l} = \kappa_l \delta_l(e^{\gamma}u_h).
\]
To obtain an explicit representation of the flux term, we approximate the diffusion and convection coefficients by constants \(\alpha_l\) and \(\bm\beta_l\), respectively, along each edge \(e_l=(i,j)\). The corresponding edgewise potential \(\gamma_l\) is defined by
\[
\gamma_l(\bm x) = \alpha_l^{-1}\bm \beta_l\cdot(\bm x- \bm x_i),
\]
so that $\nabla\gamma_l = \alpha_l^{-1}\bm \beta_l$ and $\gamma_l(\bm x_i)=0$. Introducing the Bernoulli function
\[
    \mathbf{B}(s) = \begin{cases}
        \displaystyle \frac{s}{e^s-1},&s\neq 0, \\
        1, & s=0,
    \end{cases}
\]
the edge flux can be expressed as
\begin{equation}
    \kappa_l \delta_l(e^{\gamma_l}u_h) = \alpha_l [\mathbf{B}(-p_l)u_h(\bm x_j)-\mathbf{B}(p_l)u_h(\bm x_i)],\label{eqn: edge_flux_expression}
\end{equation}
where $p_l = \alpha_l^{-1}\bm \beta_l \cdot \bm \tau_l$. Motivated by this representation, we define the Bernoulli-weighted incidence matrix \(G_b\in\mathbb R^{n_e\times n_v}\) by
\begin{align*}
    (G_b)_{lk}=\begin{cases}
        \mathbf{B}(-p_l), & \mathrm{if} \; k=j, \\
        -\mathbf{B}(p_l), & \mathrm{if} \; k = i, \\
        0, & \mathrm{otherwise},
    \end{cases}\
\end{align*}
for each oriented edge \(e_l = (i,j)\).
This modification incorporates the exponential fitting of the EAFE scheme into the discrete gradient operator.

Suppose that \(J_T(u_h)\) is constant on \(T\). Then, there exists a linear function \(\psi_T\) such that \(J_T(u_h)=\nabla\psi_T\). Taking the edge degrees of freedom and using \eqref{eqn: edge_flux_expression}, we obtain
\[
G\bm\psi_T = D_\alpha G_b\bm u_h,
\]
where \(\bm{\psi}_T\) denotes the vector of nodal values of \(\psi_T\) and 
\[
D_\alpha = \mathrm{diag}\left(\alpha_1,\dots,\alpha_{n_e}\right).
\]
Consequently, using the matrix identity \eqref{g-identity}, we have
\begin{align*}
B_T^\mathcal{N}(u_h,v_h) :=& \bm v_h^\top G^\top M D_\alpha  G_b \bm u_h\\ 
=& \bm v_h^\top G^\top M G\bm\psi_T =\bm v_h^\top G^\top W G\bm\psi_T = \bm v_h^\top G^\top W D_\alpha G_b\bm u_h.
\end{align*}
We, therefore, obtain the local EAFE bilinear form,
\begin{equation}
\label{eafe_bilinear_form_graph}
    B_T^h(u_h,v_h) = \bm v_h^\top G^\top WD_\alpha G_b\bm u_h.
\end{equation}
In the case where the flux $J_T(u_h)$ is not constant, $B_T^h$ only approximates $B_T^\mathcal{N}$.

\begin{remark}
    The graph-based formulation \eqref{eafe_bilinear_form_graph} provides a simple implementation of the EAFE scheme. Starting from the edge weights $W$ associated with the Poisson stiffness matrix, the convection-diffusion bilinear form can be assembled by approximating the coefficients as constants along each edge and incorporating the Bernoulli weights into the incidence matrix. In particular, the matrix representation \(G^\top W D_\alpha G_b\) makes explicit how the EAFE scheme approximates the flux $J(u_h)$ through an exponentially-fitted, edge-based discrete gradient. Moreover, since this formulation is constructed from barycentric coordinates and incidence relations on \(d\)-dimensional simplices, it extends the classical EAFE scheme beyond two and three dimensions while preserving its underlying algebraic structure.
\end{remark}

The extension to convex polytopal elements in arbitrary dimensions requires an additional consideration: unlike simplices, their vertex-edge graphs are not necessarily complete. Nevertheless, the incidence-based formulation \eqref{eafe_bilinear_form_graph} can be extended using generalized barycentric coordinates and complementary edges. In the following subsection, we introduce these complementary edges to construct a complete graph and extend the formulation to general convex polytopes.

\subsection{Extension to Convex Polytopal Meshes}\label{subsec: extension_polytopal}

Let $\mathcal{K}_h$ be a conforming mesh of $\Omega$ consisting of convex polytopal elements. For each $K \in \mathcal{K}_h$, let $\{\bm x_i\}_{i=1}^{n_v}\subset\mathbb{R}^d$ denote its vertices and $\{e_l\}_{l=1}^{n_e}$ its edges, where $n_v$ and $n_e$ are the numbers of vertices and edges of $K$, respectively. We orient each edge from the lower-indexed vertex to the higher-indexed vertex, as in \eqref{eq:tangent_vector}.
Let $\{\lambda_i\}_{i=1}^{n_v}$ denote the generalized barycentric coordinates on $K$, satisfying
\[
\lambda_i(\bm x_j) = \delta_{ij},\qquad \sum_{i=1}^{n_v}\lambda_i=1,\qquad \sum_{i=1}^{n_v}\lambda_i\bm x_i = \bm x.
\]
We define the corresponding conforming finite-element space by
\[
    V_h = \{v_h \in H_0^1(\Omega) : \left.v_h\right|_K \in \mathcal{W}(K), \ \forall K\in \mathcal{K}_h\},
\]
where
\[
\mathcal{W}(K)=\mathrm{span}\left\{\lambda_i\right\}_{i=1}^{n_v}.
\]

To extend the matrix identity \eqref{g-identity} to polytopal elements, let \(K^*\) denote the complete graph formed by the vertices of \(K\), with edge set \(\mathcal E_K^*\) containing
\[
n_e^* = \frac{n_v(n_v-1)}{2}
\]
edges.
Since the original edge set \(\mathcal E_K\) is a subset of \(\mathcal E_K^*\), we define the set of complementary edges by
\[
\mathcal{E}_K^c = \mathcal{E}_K^*\setminus\mathcal{E}_K.
\]
These complementary edges are auxiliary connections introduced for the local algebraic construction and do not modify the underlying polytopal mesh.
We define the vertex-edge incidence matrix \(G\in\mathbb R^{n_e^*\times n_v}\) on the complete graph \(K^*\) by
\[
    G_{lk} = \begin{cases}
        1, & \mathrm{if} \; k=j, \\
        -1, & \mathrm{if} \; k = i, \\
        0, & \mathrm{otherwise},
    \end{cases}
\]
for each oriented edge \(e_l=(i,j)\in\mathcal E_K^*\).

Associated with each edge of the complete graph, including complementary edges, we define the generalized Whitney edge function
\begin{equation}
    \bm \phi_l = \lambda_i\nabla \lambda_j - \lambda_j \nabla \lambda_i,\qquad i<j.\label{eqn: generalized_Whitney_edge}
\end{equation}
Unlike the simplex case, these generalized Whitney edge functions do not necessarily form a canonical basis with respect to the edge degrees of freedom.
However, since the generalized barycentric coordinates satisfy the partition-of-unity property, the same argument as in the simplex case in \eqref{eq:edge_relations} gives
\begin{equation}
\nabla \lambda_i = \sum_{l=1}^{n_e^*}G_{li}\bm\phi_l.\label{eqn: incidence-based_gradient_polytopes}  
\end{equation}
Thus, introducing complementary edges preserves the incidence-based representation of the nodal basis gradients. This allows us to extend the matrix identity \eqref{g-identity} to convex polytopal elements.
We first illustrate this construction with an example before establishing the general result.

\begin{figure}[ht!]
\centering
\begin{tikzpicture}[scale=2.5]

\coordinate (A1) at (0,0);
\coordinate (A2) at (1,0);
\coordinate (A4) at (1,1);
\coordinate (A3) at (0,1);

\draw[->, thick] (A1) -- (A2);
\node[below] at (0.5,0) {$e_1$};

\draw[->, thick] (A1) -- (A3);
\node[left] at (0,0.5) {$e_2$};

\draw[->, thick] (A2) -- (A4);
\node[right] at (1,0.5) {$e_3$};

\draw[->, thick] (A3) -- (A4);
\node[above] at (0.5,1) {$e_4$};

\foreach \v in {A1,A2,A3,A4}
    \fill (\v) circle (0.02);

\node[below left]  at (A1) {$1$};
\node[below right] at (A2) {$2$};
\node[above left]  at (A3) {$3$};
\node[above right] at (A4) {$4$};

\coordinate (B1) at (2.5,0);
\coordinate (B2) at (3.5,0);
\coordinate (B4) at (3.5,1);
\coordinate (B3) at (2.5,1);

\draw[->, thick] (B1) -- (B2);
\node[below] at (3,0) {$e_1$};

\draw[->, thick] (B1) -- (B3);
\node[left] at (2.5,0.5) {$e_2$};

\draw[->, thick] (B2) -- (B4);
\node[right] at (3.5,0.5) {$e_3$};

\draw[->, thick] (B3) -- (B4);
\node[above] at (3,1) {$e_4$};

% Diagonals

\draw[->, dashed] (B1) -- (B4);
\node at (2.9,0.25) {$e_5$};

\draw[->, dashed] (B2) -- (B3);
\node[left] at (3,0.75) {$e_6$};

\foreach \v in {B1,B2,B3,B4}
    \fill (\v) circle (0.02);

\node[below left]  at (B1) {$1$};
\node[below right] at (B2) {$2$};
\node[above left]  at (B3) {$3$};
\node[above right] at (B4) {$4$};

\node at (0.5,-0.3) {$K$};
\node at (3,-0.3) {$K^*$};
\end{tikzpicture}
\caption{
Original quadrilateral $K$ (left) and its complete graph $K^*$ (right). Dashed lines represent complementary edges.
}
\label{diagram}
\end{figure}
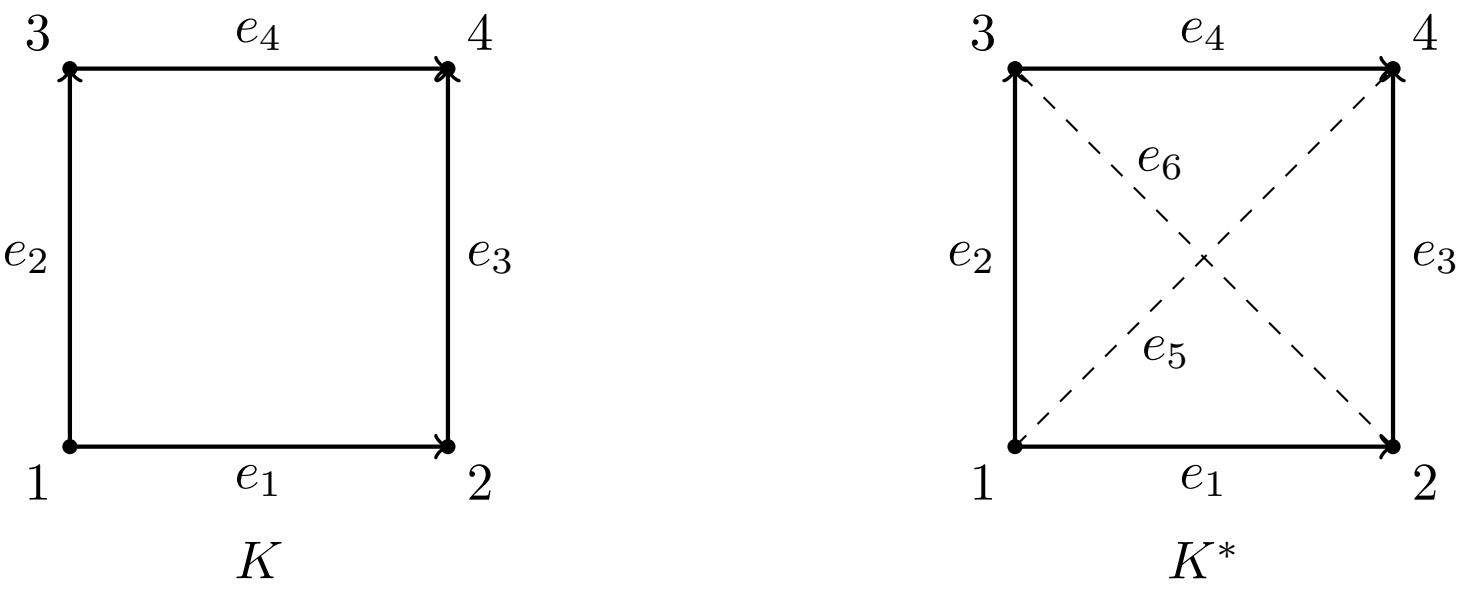

\textit{\textbf{Example.}}
Consider a quadrilateral element $K$ with four vertices and four edges, as illustrated in \Cref{diagram}. Its vertex-edge graph is not complete, since $n_v=4$ and $n_e=4$, whereas the complete graph contains $n_e^* = 4(4-1)/2=6$ edges. We therefore introduce two complementary edges, $e_5$ and $e_6$, connecting the two pairs of nonadjacent vertices.

Following the edge orientation convention, we assign:
\begin{align*}
    &\text{Original:}\quad e_1 = (1, 2),\quad e_2 = (1,3),\quad e_3 = (2, 4),\quad e_4 = ( 3, 4),\\
    &\text{Complementary:}\quad e_5 = (1,4),\quad e_6 = (2,3).
\end{align*}
The corresponding incidence matrix of the complete graph is
\[
G=
\begin{bmatrix}
-1& 1& 0& 0\\
-1& 0& 1& 0\\
0&-1& 0& 1\\
0& 0&-1& 1\\[-1.5ex]
\hdotsfor{4}\\[-0.5ex]
-1& 0& 0& 1\\
0&-1& 1& 0
\end{bmatrix}.
\]
The first four rows correspond to the original edges of \(K\), while the last two rows correspond to the complementary edges.

For the unit square \(K=[0,1]^2\), the generalized barycentric coordinates are given by the bilinear Lagrange basis functions
\begin{align*}
\lambda_1(x,y) = (1-x)(1-y),\qquad&\lambda_2(x,y) = x(1-y),\\
\lambda_3(x,y) = (1-x)y,\qquad&\lambda_4(x,y) = xy,
\end{align*}
satisfying the partition-of-unity property.
The corresponding generalized Whitney edge functions $\bm\phi_l$ are constructed as in \eqref{eqn: generalized_Whitney_edge}. In particular, the two complementary edges give rise to the additional edge functions
\[
\bm\phi_5 = \lambda_1\nabla\lambda_4 - \lambda_4\nabla\lambda_1,\qquad
\bm\phi_6 = \lambda_2\nabla\lambda_3 - \lambda_3\nabla\lambda_2.
\]
Thus, the example incidence matrix implies
\[
\nabla\lambda_i = \sum_{l=1}^6G_{li}\bm\phi_{l},\qquad i=1,2,3,4,
\]
so the complementary edges provide the additional terms needed to recover the incidence-based representation of the nodal basis gradients. 

This construction extends beyond rectangular elements. For general convex polytopes, generalized barycentric coordinates, such as Wachspress coordinates \cite{wachspress3d,wachspress2011barycentric}, can be used to construct the corresponding nodal and generalized Whitney edge functions. The same incidence-based gradient identity then follows from the partition-of-unity property on the complete graph.

\subsection{Matrix Identity on Convex Polytopes}
We now establish the polytopal counterpart of the matrix identity \eqref{g-identity} using the complementary edges introduced in \Cref{subsec: extension_polytopal}.
\begin{lemma} \label{m and w}
Let \(G\in\mathbb R^{n_e^*\times n_v}\) be the incidence matrix of the complete graph \(K^*\), and let \(M\in\mathbb R^{n_e^*\times n_e^*}\) be the generalized Whitney edge mass matrix defined by
\[
M_{lm} = \int_K\bm\phi_l\cdot\bm\phi_m\dx.
\]
Then,
    \begin{equation}\label{eqn: matrix_identity_polytopes}
        G^\top M G = G^\top WG,
    \end{equation}
where \(W\in\mathbb R^{n_e^*\times n_e^*}\) is the diagonal matrix with entries
    \[
        W_{ll} = -\int_K \nabla \lambda_i \cdot \nabla \lambda_j\dx,\qquad e_l = (i,j).
    \]
\end{lemma}
\begin{proof}
We establish the identity by comparing the entries of \(G^\top MG\) and \(G^\top WG\).
For \(i\neq j\), the complete graph \(K^*\) contains exactly one edge connecting vertices \(i\) and \(j\). Since \(G_{li}G_{lj}=-1\) for this edge and vanishes otherwise, we
have
\[
(G^\top WG)_{ij}=\sum_{l=1}^{n_e^*}G_{li}W_{ll}G_{lj}
=-W_{ll}
=\int_K\nabla\lambda_i\cdot\nabla\lambda_j\dx.
\]
For \(i=j\), using the partition-of-unity property of the generalized barycentric coordinates, we obtain
\[
(G^\top WG)_{ii}=\sum_{l=1}^{n_e^*}G_{li}^2W_{ll}=-\sum_{\substack{k=1\\k\neq i}}^{n_v}
\int_K\nabla\lambda_i\cdot\nabla\lambda_k\dx=\int_K|\nabla\lambda_i|^2\dx.
\]
Thus, for all $i,j$,
\[
(G^\top WG)_{ij}=\int_K\nabla\lambda_i\cdot\nabla\lambda_j\dx.
\]

On the other hand, using the incidence-based gradient representation \eqref{eqn: incidence-based_gradient_polytopes}, we have
\begin{align*}
(G^\top MG)_{ij}
&=\sum_{l=1}^{n_e^*}\sum_{m=1}^{n_e^*}G_{li}M_{lm}G_{mj}\\
&=\int_K
\left(\sum_{l=1}^{n_e^*}G_{li}\boldsymbol\phi_l\right)
\cdot
\left(\sum_{m=1}^{n_e^*}G_{mj}\boldsymbol\phi_m\right)\dx\\
&=\int_K\nabla\lambda_i\cdot\nabla\lambda_j\dx.
\end{align*}
Therefore,
\[
G^\top M G = G^\top W G.
\]

\end{proof}
The identity extends the diagonal representation of the Poisson stiffness matrix from simplices to convex polytopal elements, generalizing \eqref{g-identity}. It provides the algebraic foundation for the edge-averaged discretization developed in the following section.

Recall from the example in Figure~\ref{diagram} where we introduce two new edges with corresponding basis functions, $\bm \phi_5$ and $\bm \phi_6$. This allows us to establish the necessary relationship between the edge basis functions and the barycentric coordinates in order to ensure the matrix identity in Lemma~\ref{m and w}. However, in the implementation, since we replace the mass matrix with the weight matrix, there is no reason to ever construct the edge basis explicitly, and we can instead rely entirely on the barycentric coordinates in the original graph.

\section{Incidence-Based Edge-Averaged Method}\label{section:scheme}

We next extend the EAFE scheme to convex polytopal meshes using the incidence-based framework developed in \Cref{sec:graph}. By combining the Bernoulli-weighted incidence matrix with the diagonal representation of the Poisson stiffness matrix, we construct the incidence-based edge-averaged method (IBEAM) for convex polytopes in arbitrary dimensions. We then establish its monotonicity under suitable conditions on the edge weights.

Following the EAFE construction in \Cref{subsec: EAFE}, we approximate the flux $J(u_h)$ locally by $J_K(u_h)$ and consider the corresponding bilinear form
\[
B_K^\mathcal{N}(u_h,v_h)=\int_K J_K(u_h)\cdot\nabla v_h\dx.
\]
For each oriented edge \(e_l=(i,j)\) of the complete graph \(K^*\), including complementary edges, we approximate the tangential flux by a constant value along \(e_l\). Using the exponential-fitting construction in \eqref{eqn: exponential_fitting}, we impose the weighted edge degrees of freedom
\[
\mathrm{dof}_l(\alpha^{-1}e^\gamma J_K(u_h))=\mathrm{dof}_l(\alpha^{-1}e^\gamma J(u_h))=\delta_l(e^\gamma u_h).
\]
Since $J_K(u_h)\cdot\bm\tau_l$ is constant along $e_l$, it follows that
\[
\mathrm{dof}_l(J_K(u_h))=\left.(J_K(u_h)\cdot\bm\tau_l)\right|_{e_l} = \kappa_l\delta_l(e^{\gamma} u_h).
\]
Using the edgewise constant-coefficient approximations \(\alpha_l\) and \(\bm{\beta}_l\), together with \eqref{eqn: edge_flux_expression}, we obtain
\begin{equation}
    \mathrm{dof}_l(J_K(u_h))=\alpha_l [\mathbf{B}(-p_l)u_h(\bm x_j)-\mathbf{B}(p_l)u_h(\bm x_i)] = (D_\alpha G_b\bm u_h)_l,\label{eqn: edge_flux_vector_form}
\end{equation}
where $p_l = \alpha_l^{-1}\bm\beta_l\cdot\bm\tau_l$. Thus, the exponentially-fitted edge fluxes are represented by the vector \(D_\alpha G_b\bm u_h\).

We now establish an identity for the local bilinear form under the assumption of an elementwise constant approximated flux.
\begin{lemma}\label{lemma: constant_flux_identity}
Suppose that \(J_K(u_h)\) is constant on \(K\) and that its edge degrees of freedom are given by the exponentially-fitted fluxes in \eqref{eqn: edge_flux_vector_form}. Then,
\[
B_K^\mathcal{N}(u_h,v_h) = \bm v_h^\top G^\top W D_\alpha G_b\bm u_h.
\]
\end{lemma}
\begin{proof}
    Since \(J_K(u_h)\) is constant on \(K\), there exists an affine function \(\psi_K\in\mathcal W(K)\) such that
    \[
    J_K(u_h) = \nabla\psi_K.
    \]
    For each oriented edge \(e_l=(i,j)\), applying the edge degree of freedom gives
    \[
    \mathrm{dof}_l(J_K(u_h))=\frac{1}{|e_l|}\int_{e_l}\nabla\psi_K\cdot\bm\tau_l\ds = \psi_K(\bm x_j) - \psi_K(\bm x_i) = (G\bm\psi_K)_l,
    \]
    where $\bm\psi_K$ denotes the vector of nodal values of $\psi_K$. Combining this relation with \eqref{eqn: edge_flux_vector_form}, we obtain
    \[
    G\bm\psi_K = D_\alpha G_b\bm u_h.
    \]
    Consequently, using the incidence-based gradient representation and \Cref{m and w}, we have
    \begin{align*}
        B_K^\mathcal{N}(u_h,v_h)&= \int_K\nabla\psi_K\cdot\nabla v_h\dx\\
        &=\bm v_h^\top G^\top M G \bm\psi_K\\
        &=\bm v_h^\top G^\top W G \bm\psi_K\\
        &=\bm v_h^\top G^\top WD_\alpha G_b \bm u_h,
    \end{align*}
    which completes the proof.
\end{proof}

Again, for a general approximated flux \(J_K(u_h)\), motivated by the constant-flux case, we define the local IBEAM bilinear form by
\begin{equation}\label{ibeam_bi}
    B_K^h(u_h,v_h): = \bm v_h^\top G^\top W D_\alpha G_b\bm u_h.
\end{equation}
The proposed IBEAM scheme is then to find $u_h\in V_h$ such that
\[
\sum_{K\in\mathcal{K}_h} B_K^h(u_h,v_h) = F(v_h),\qquad \forall v_h\in V_h.
\]
The corresponding local stiffness matrix is
\[
A_K = G^\top W D_\alpha G_b,
\]
and the global stiffness matrix \(A\) is assembled from these local contributions in the standard finite-element manner.

\subsection{A Sufficient Condition for Monotonicity Property}\label{section:monotone}

Finally, we establish a sufficient condition for the monotonicity of the proposed IBEAM scheme. Specifically, we show that the resulting stiffness matrix is a nonsingular $M$-matrix under suitable positivity and connectivity conditions.

\begin{lemma}
    Suppose that \(\omega_l^K>0\) for every edge \(e_l\) of the complete graph \(K^*\) associated with each \(K\in\mathcal K_h\), and that \(\alpha_l>0\). Then, the global stiffness matrix \(A\), restricted to the interior degrees of freedom, is a nonsingular $M$-matrix.
\end{lemma}

\begin{proof}
    The global stiffness matrix is assembled from local contributions of the form \(G^\top WD_\alpha G_b\). For an oriented edge \(e_l=(i,j)\), the corresponding contribution, restricted to the two vertices \(i\) and \(j\), is
\[
A_K^{(l)}
=
\omega_l^K\alpha_l
\begin{bmatrix}-1\\1\end{bmatrix}
\begin{bmatrix}
-\mathbf B(p_l)&\mathbf B(-p_l)
\end{bmatrix}=
\omega_l^K\alpha_l
\begin{bmatrix}
\mathbf B(p_l)&-\mathbf B(-p_l)\\
-\mathbf B(p_l)&\mathbf B(-p_l)
\end{bmatrix}.
\]
Since \(\mathbf B(s)>0\) for all \(s\in\mathbb R\), \(\omega_l^K>0\), and \(\alpha_l>0\) by assumption, the diagonal entries of each edge contribution are positive, while its off-diagonal entries are negative. Therefore, the assembled global stiffness matrix has nonpositive off-diagonal entries and is a $Z$-matrix.

Summing the columns of each local edge contribution gives
\[
\begin{bmatrix}1&1\end{bmatrix}
\begin{bmatrix}
\mathbf B(p_l)&-\mathbf B(-p_l)\\
-\mathbf B(p_l)&\mathbf B(-p_l)
\end{bmatrix}
=
\begin{bmatrix}0&0\end{bmatrix}.
\]
Thus, the fully assembled stiffness matrix $\bar{A}$, before imposing Dirichlet boundary conditions, satisfies
\[
\bm 1^\top \bar{A}
=
\bm 0^\top.
\]
Since \(\bar{A}\) is a $Z$-matrix with positive diagonal entries, it is weakly column-diagonally dominant.
After imposing homogeneous Dirichlet boundary conditions, the reduced matrix \(A\) remains a weakly column-diagonally dominant \(Z\)-matrix, with strict dominance in columns corresponding to interior vertices adjacent to the boundary. Since every connected component of the interior vertex graph is connected to the Dirichlet boundary, \(A\) is a nonsingular \(M\)-matrix by the irreducible diagonal dominance criterion applied to each component.
\end{proof}

\begin{remark}
    The nonsingular \(M\)-matrix property implies that \(A^{-1}\geq0\) entrywise, ensuring the existence and uniqueness of the discrete IBEAM solution and the discrete maximum principle. However, this monotonicity result does not by itself guarantee convergence of the proposed scheme. On simplicial meshes, IBEAM extends the classical EAFE scheme from two and three dimensions to arbitrary dimensions. Since the underlying finite-element structure remains unchanged, the first-order convergence analysis of EAFE can be extended to arbitrary-dimensional simplices under analogous regularity and mesh assumptions. We therefore omit the proof. For general convex polytopal meshes, however, the same argument cannot be applied directly, since the generalized barycentric basis functions are not necessarily piecewise linear and may have nonzero second derivatives. Nevertheless, the numerical experiments in the following section demonstrate first-order convergence in the \(H^1\) semi-norm on various polytopal meshes.
\end{remark}

\section{Numerical Experiments}\label{section:numerical}

To investigate the efficacy of the IBEAM method, we perform two types of numerical experiments in two and three dimensions. The first is a convection-dominated problem with a boundary layer and no closed-form solution, which we use to assess qualitatively whether the scheme avoids the spurious oscillations that afflict the standard Galerkin method. The second set of tests is a convergence study with a manufactured solution, in which we measure the $L^2$ norm and $H^1$ semi-norm errors under mesh refinement. In the latter set of test problems, we simulate on square, quadrilateral, and polygonal meshes in two dimensions and on tetrahedral and hexahedral meshes in three dimensions. We again note that polygonal and hexahedral meshes have not been previously treated by EAFE-type methods \cite{safe} (nor hexahedral for EAVE-type schemes \cite{eave}), and therefore we highlight the performance of IBEAM on these meshes.  We include tetrahedral meshes as a validation that pre-existing results would remain consistent.  

All computations were performed in MATLAB using Wachspress coordinates \cite{wachspress3d,wachspress2011barycentric} for the barycentric coordinates and a five-point Gaussian quadrature rule. A direct solver was used to solve the resulting linear systems.  Simulations were run on the CPU of a System76 Pangolin laptop with an AMD Ryzen 7 6800U processor (8 cores, 16 threads, up to 4.7~GHz) and 32~GB of LPDDR5 memory.

\subsection{Two-Dimensional Results}\label{section:2d_results}
\begin{figure}[h!]
    \centering
    % Top row
    \begin{subfigure}{0.3\linewidth}
        \centering
        \includegraphics[scale=0.32]{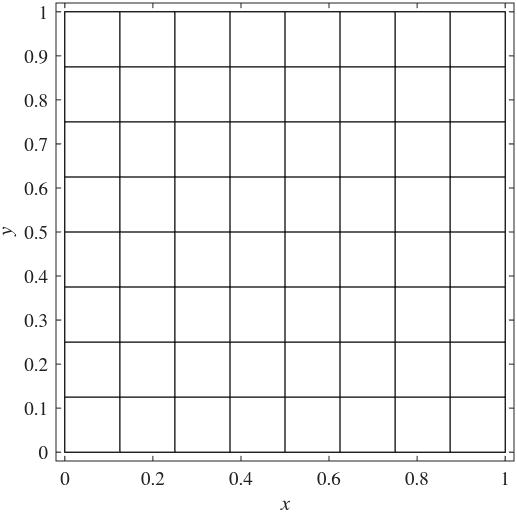}
        \caption{Square mesh.}
        \label{fig:square_mesh}
    \end{subfigure}
    \hskip 10pt
    \begin{subfigure}{0.3\linewidth}
        \centering
        \includegraphics[scale=0.32]{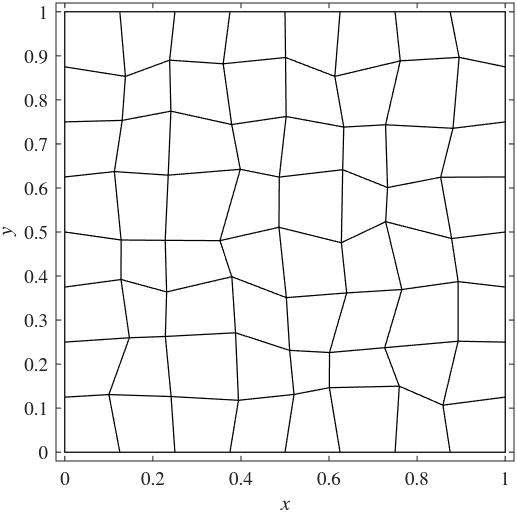}
        \caption{Quadrilateral mesh.}
        \label{fig:quad_mesh}
    \end{subfigure}
    \hskip 10pt
    \begin{subfigure}{0.3\linewidth}
        \centering
        \includegraphics[scale=0.32]{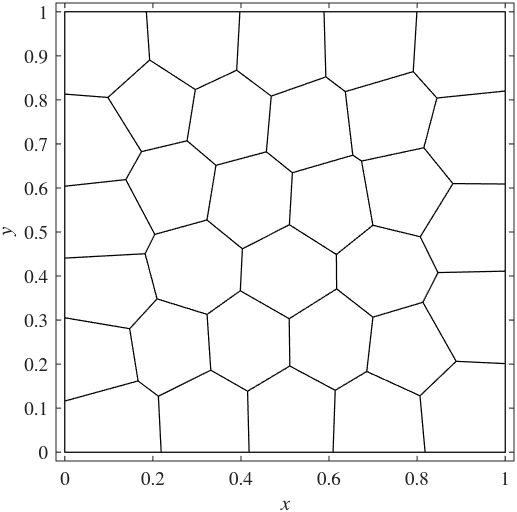}
        \caption{Polygonal mesh.}
        \label{fig:poly_mesh}
    \end{subfigure}    
    \caption{Two-dimensional meshes used in numerical experiments.}
    \label{fig:meshes}
\end{figure}

For the two-dimensional results, we consider three mesh types: square, quadrilateral, and general polygonal. We choose the square mesh for ease of construction while still demonstrating the method's effectiveness for non-simplicial grids. Quadrilaterals are chosen to see how the method is affected when presented with an unstructured mesh, and convex polygons reinforce the schemes performance on more general meshes. We present a visual representation of each of these grids in \Cref{fig:meshes}.
\subsubsection{Boundary Layer Example}\label{section:2d_boundary_layer} Our first test is a convection-dominated convection-diffusion problem with $\alpha = 10^{-3}$ and $10^{-6}$, $\bm\beta = [1\;2]^\top$, and $f = 1$. Problems of this type arise in many applications, including heat and mass transport, pollutant dispersion, and semiconductor device modeling, where diffusion is weak relative to convection and the solution develops thin boundary layers.
The relative strength of convection and diffusion is measured by the mesh P\'eclet number,
\[
  \mathrm{Pe}_h = \frac{|\bm\beta|h}{2\alpha},
\]
which compares the amount of convective transport across one element to the amount of diffusive transport over the same distance. When $\mathrm{Pe}_h \le 1$, diffusion dominates on the mesh scale, and a standard Galerkin method performs well. When $\mathrm{Pe}_h > 1$, convection dominates at the mesh scale, the stiffness matrix may lose its $M$-matrix property, and Galerkin solutions develop spurious oscillations \cite{survey}. Equivalently, $\mathrm{Pe}_h > 1$ means the mesh is too coarse to resolve the boundary layer, whose width is of order $\alpha/|\bm\beta|$ \cite{roos2008robust}. Large P\'eclet numbers suggest a demanding test of the monotonicity properties of the IBEAM system established in \Cref{section:monotone}.

For the following experiments, we compute a solution using IBEAM on a square mesh and the standard Galerkin method on a uniform triangular mesh, both with a mesh size of $h = 2^{-5}$. For the present parameters, $\mathrm{Pe}_h = \sqrt{5}\cdot 2^{-5}/(2\cdot 10^{-3}) \approx 35$ for $\alpha = 10^{-3}$ and $35,000$ for $\alpha = 10^{-6}$, so the problem is well inside the convection-dominated regime. 

\begin{figure}[h!]
    \centering
    \begin{subfigure}{0.9\linewidth}
        \centering
        \includegraphics[scale=0.32]{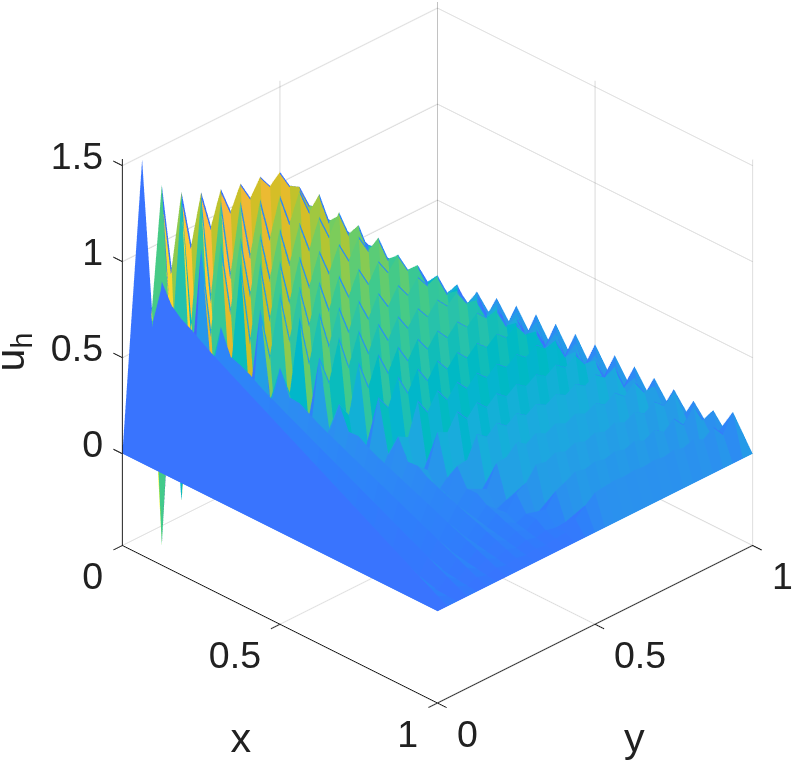}
        \caption{FEM solution, $\alpha=10^{-3}$.}
        \label{ex1fem}
    \end{subfigure}
    \hskip 30pt
    \begin{subfigure}{0.45\linewidth}
        \centering
        \includegraphics[scale=0.22]{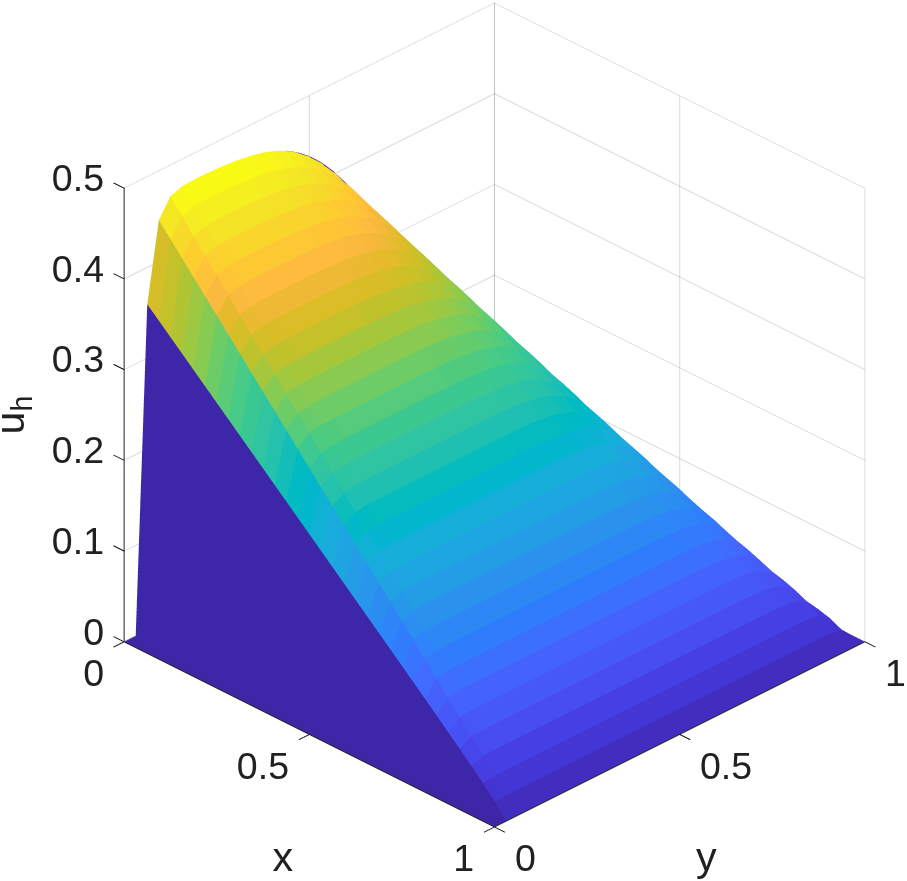}
        \caption{IBEAM solution, $\alpha=10^{-3}$.}
        \label{ex2graph}
    \end{subfigure}
    \hskip 30pt
    \begin{subfigure}{0.45\linewidth}
        \centering
        \includegraphics[scale=0.22]{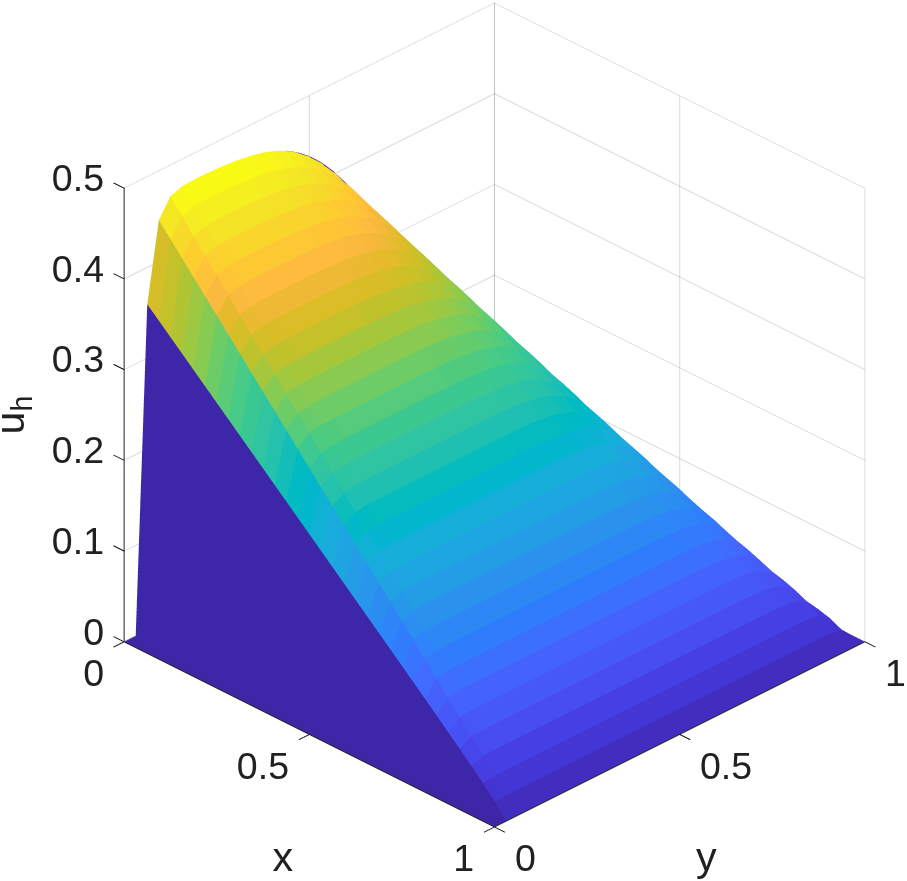}
        \caption{IBEAM solution, $\alpha=10^{-6}$.}
        \label{ex3graph}
    \end{subfigure}
    \vskip -10pt
    \caption{Comparison of FEM and IBEAM  solutions for $\alpha = 10^{-3}$, $\bm \beta = [1 \; 2]^\top$,
    $f=1$, and $h=1/32$, where the IBEAM solution refers to the solution
    utilizing local stiffness matrix \eqref{ibeam_bi}. Unlike the FEM solution, which exhibits
    nonphysical oscillations, the IBEAM solution
    avoids oscillatory behavior, even for significantly smaller $\alpha=10^{-6}$.}
    \label{ex1comparison}
\end{figure}

\Cref{ex2graph,ex3graph} show that the IBEAM solution is free of spurious oscillations, consistent with the known behavior of EAFE-type schemes \cite{eafe} and with the $M$-matrix and discrete maximum principle properties in \Cref{section:monotone}. By contrast, the standard Galerkin method oscillates near the boundary layer \cite{survey}, even with the larger value of $\alpha$, rendering the solution nonphysical.

\subsubsection{Convergence Study}\label{section:2d_convergence}
Next, we perform a convergence study using the manufactured solution $u(x,y) = \sin(\pi x)\sin(\pi y)$, with convection vector $\bm\beta = [0\;1]^\top$ and right-hand side chosen so that $u$ satisfies the convection-diffusion equation. The diffusion coefficient takes values $\alpha \in \{1, 0.1, 0.01, 0.001\}$, and the mesh size varies from $h = 2^{-2}$ to $h = 2^{-8}$. For each mesh family discussed above, we report the $L^2$ norm and $H^1$ semi-norm error at each refinement level. We compare the observed $H^1$ semi-norm rates with those known for EAFE: an $O(h)$ estimate is proven in \cite{eafe}.

\newlength{\errw}
\setlength{\errw}{1.6cm}
\newlength{\ordw}
\setlength{\ordw}{0.7cm}

\begin{figure}[h!]
    \centering
    \begin{subfigure}{0.45\linewidth}
        \centering
        \includegraphics[width=\linewidth]{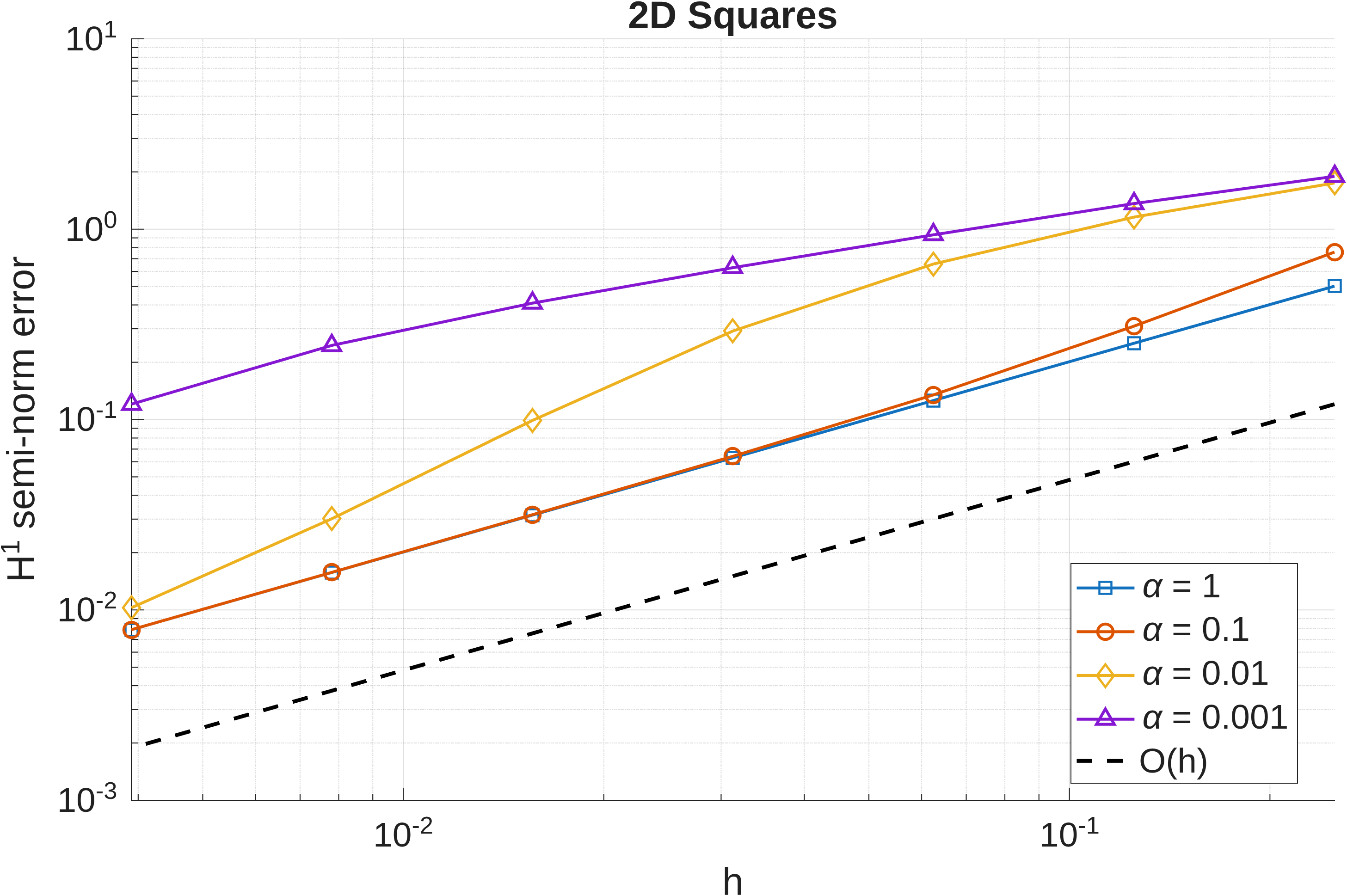}
        %\caption{$H^1$ semi-norm error.}
        \label{fig:conv_2d_square_h1}
    \end{subfigure}
    \hskip 30pt
    \begin{subfigure}{0.45\linewidth}
        \centering
        \includegraphics[width=\linewidth]{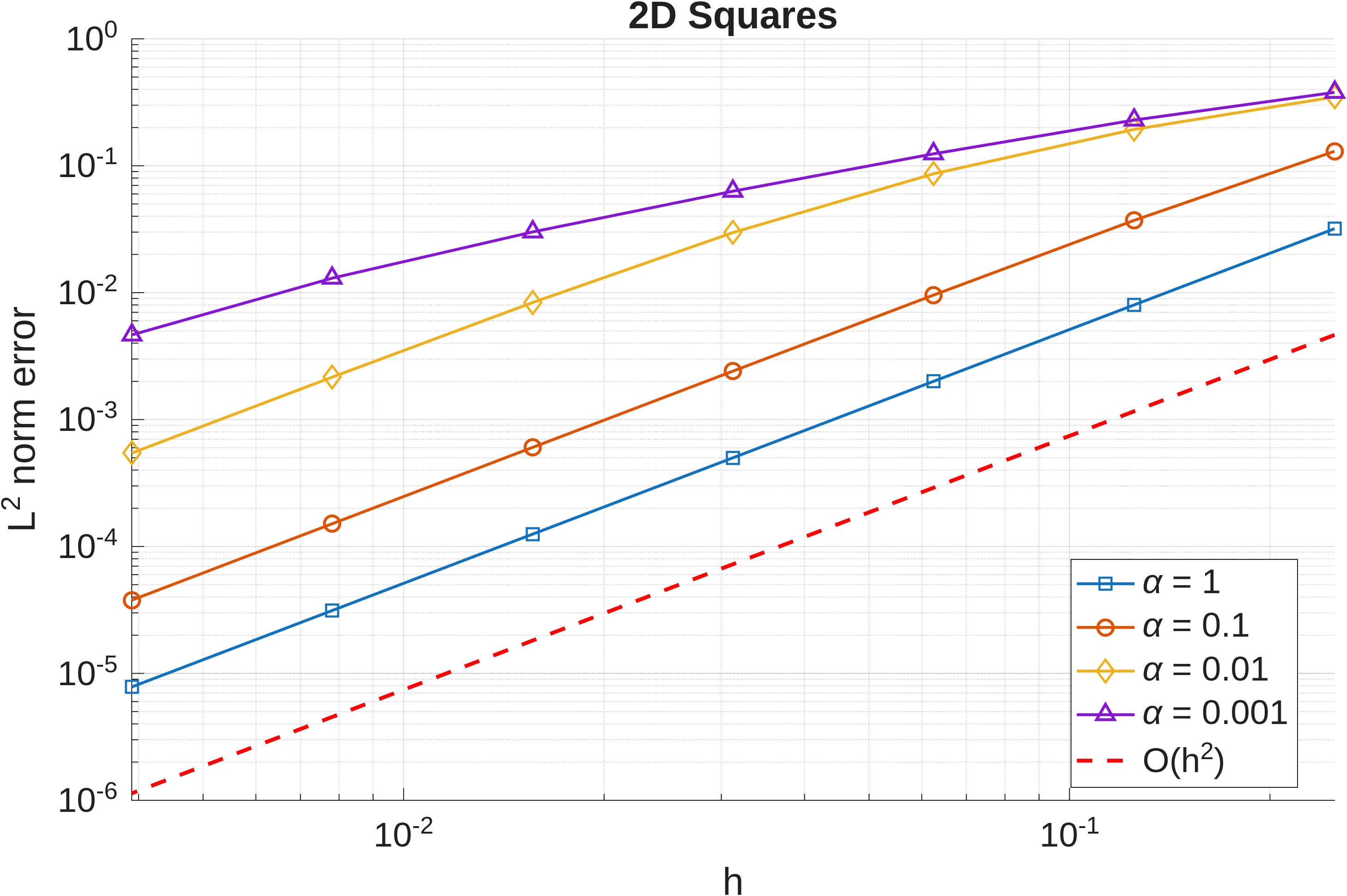}
        %\caption{$L^2$ norm error.}
        \label{fig:conv_2d_square_l2}
    \end{subfigure}
    \vskip -10pt
    \caption{IBEAM error convergence rates on square meshes for $\alpha \in \{1, 0.1, 0.01, 0.001\}$, plotted against mesh size $h$, for the $H^1$ semi-norm (left) and  $L^2$ norm (right) error.}
    \label{fig:conv_2d_square}
\end{figure}

On square meshes (\Cref{fig:conv_2d_square}), in the diffusion-dominated regime ($\alpha = 1$), the computed rates agree with the established rate \cite{eafe} of $O(h)$ in the $H^1$ semi-norm from the coarsest mesh onward, and we observe $O(h^2)$ convergence in the $L^2$ norm. As $\alpha$ decreases, the order is reduced on coarse meshes, and additional refinement is needed before the optimal rate is reached. We call this range of mesh sizes the pre-asymptotic regime. 

\begin{figure}[h!]
    \centering
    \begin{subfigure}{0.45\linewidth}
        \centering
        \includegraphics[width=\linewidth]{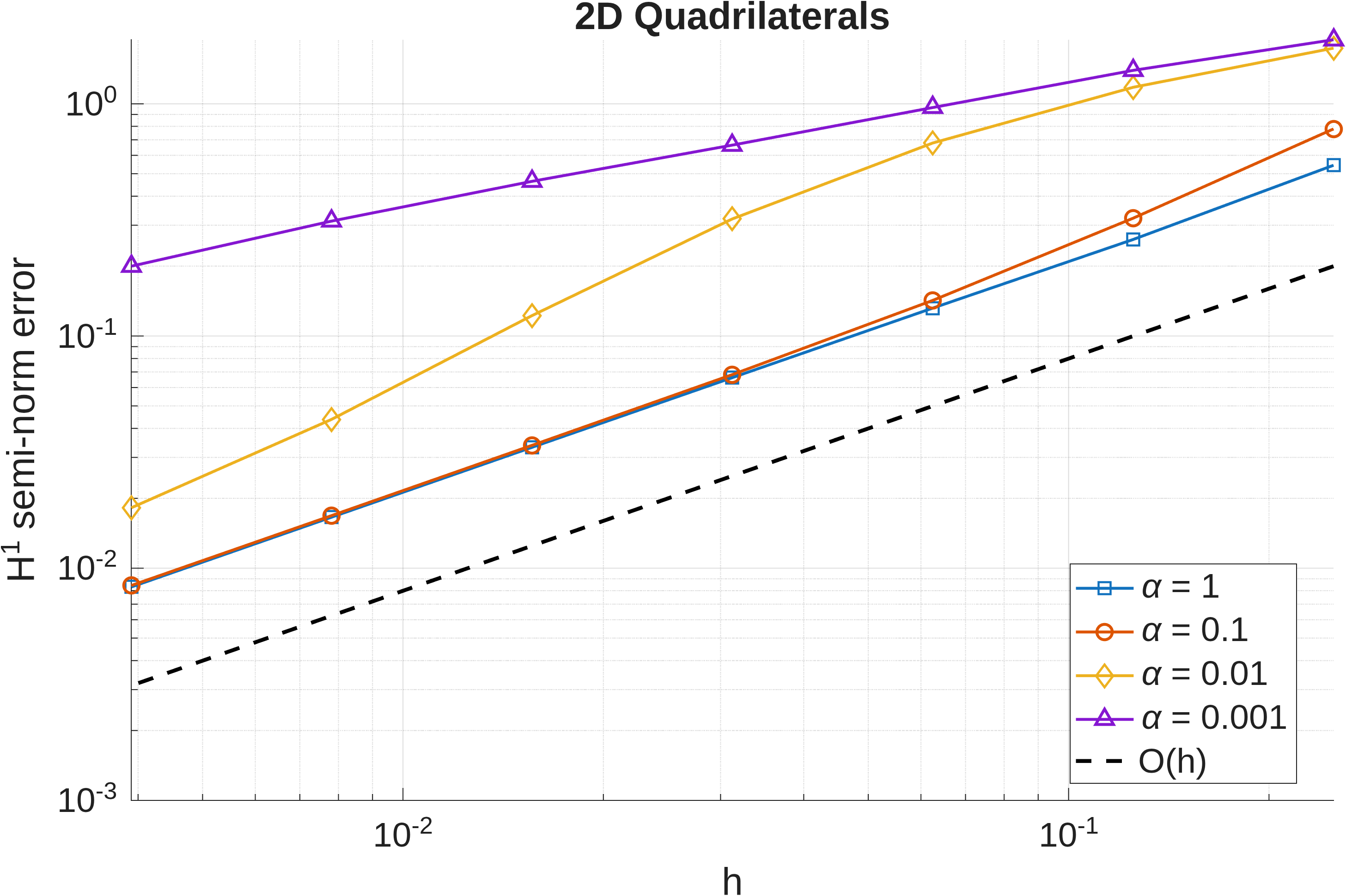}
        \label{fig:conv_2d_quad_h1}
    \end{subfigure}
    \hskip 30pt
    \begin{subfigure}{0.45\linewidth}
        \centering
        \includegraphics[width=\linewidth]{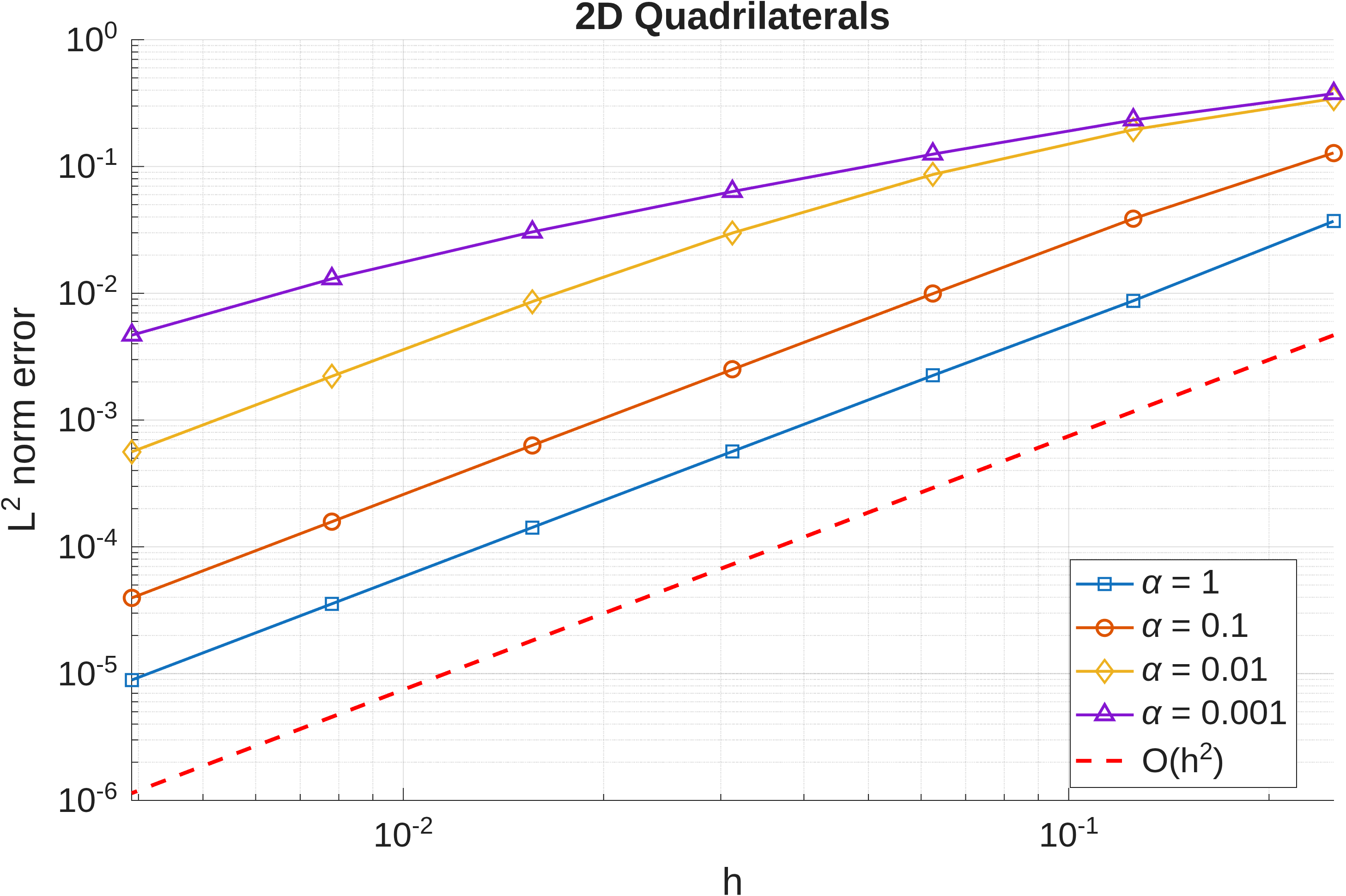}
        \label{fig:conv_2d_quad_l2}
    \end{subfigure}
    \vskip -10pt
    \caption{IBEAM error convergence rates on quadrilateral meshes for $\alpha \in \{1, 0.1, 0.01, 0.001\}$, plotted against mesh size $h$, for the $H^1$ semi-norm (left) and  $L^2$ norm (right) error.}
    \label{fig:conv_2d_quad}
\end{figure}

\begin{figure}[h!]
    \centering
    \begin{subfigure}{0.45\linewidth}
        \centering
        \includegraphics[width=\linewidth]{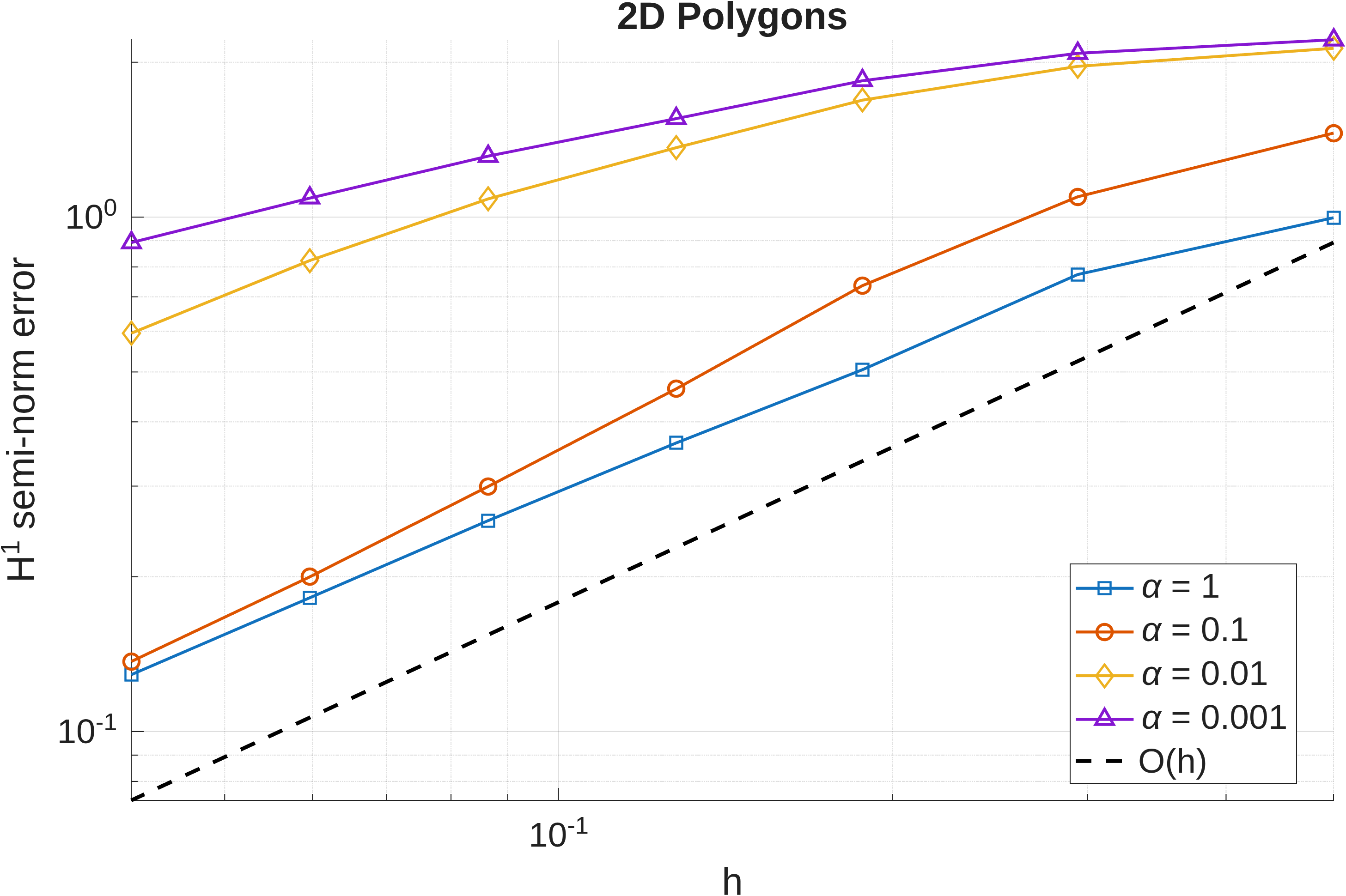}
        \label{fig:conv_2d_poly_h1}
    \end{subfigure}
    \hskip 30pt
    \begin{subfigure}{0.45\linewidth}
        \centering
        \includegraphics[width=\linewidth]{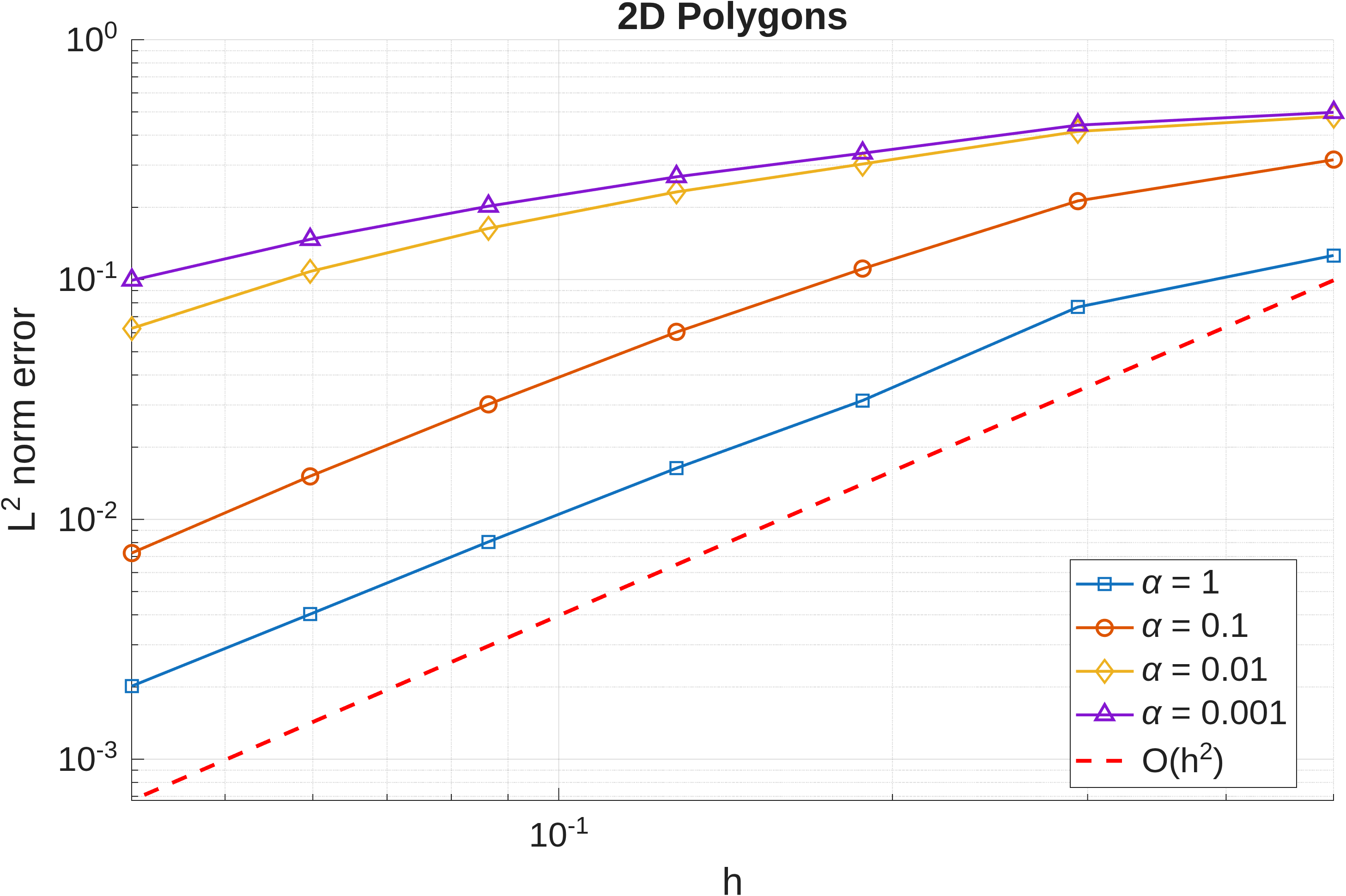}
        \label{fig:conv_2d_poly_l2}
    \end{subfigure}
    \vskip -10pt
    \caption{IBEAM error convergence rates on convex polygonal meshes for $\alpha \in \{1, 0.1, 0.01, 0.001\}$, plotted against mesh size $h$, for the $H^1$ semi-norm (left) and  $L^2$ norm (right) error.}
    \label{fig:conv_2d_poly}
\end{figure}

On quadrilateral meshes (\Cref{fig:conv_2d_quad}), the pre-asymptotic regime of convergence appears at $\alpha = 10^{-1}$, earlier than on square meshes. As the mesh becomes less structured, we expect the condition number of the element stiffness matrix to increase and the conditions required for monotonicity, as shown in \Cref{section:monotone}, to be violated. For general convex polygonal meshes (\Cref{fig:conv_2d_poly}), we approximate the mesh size $h$ by $N_{\mathrm{el}}^{-1/2}$. Here, we see that the pre-asymptotic regime is now noticeable across all diffusion values, requiring further refinement to correct it. The difference from a structured mesh is most pronounced: some elements do not satisfy the conditions established for monotonicity in \Cref{section:monotone}. In particular, large variations in edge lengths and interior angles cause stability issues.
This manifests as a degradation in quality for the $L^2$ norm convergence order, but the $O(h)$ convergence in the $H^1$ semi-norm is partially preserved. In particular, for $\alpha = 10^{-3}$, the error does not yet appear to converge at $O(h)$, but we expect the asymptotic regime to emerge with further mesh refinement.

\subsection{Three-Dimensional Results}\label{section:3d_results}
For the three-dimensional results, we consider two mesh types: tetrahedral and hexahedral. The tetrahedral mesh allows for direct comparison to a variety of schemes that function on simplicial meshes and serves as a baseline result. The hexahedral mesh is chosen for its simplicity of construction and to demonstrate the scheme's performance on non-simplicial meshes.

\subsubsection{Boundary Layer Example}

We return to the boundary layer model problem, but in three dimensions. For the following experiments, we compute a solution using IBEAM on a hexahedral mesh and the standard Galerkin method on a uniform tetrahedral mesh. The parameter regime is kept consistent: $\alpha = 10^{-3}$ and $\alpha = 10^{-6}$, $f = 1$, and $\bm \beta = [1\;2\;0]^\top$, but we reduce the mesh resolution by a factor of two, $h = 2^{-4}.$ The mesh P\'eclet numbers associated with these parameter regimes are $\mathrm{Pe}_h = \sqrt{5}\cdot 2^{-4}/(2\cdot 10^{-3}) \approx 70$ for $\alpha=10^{-3}$ and $70,000$ for $10^{-6}$, so we expect a boundary layer to form which the standard Galerkin finite-element method would fail to resolve. 

\begin{figure}[h!]
    \centering
    \begin{subfigure}{0.9\linewidth}
        \centering
        \includegraphics[scale=0.28]{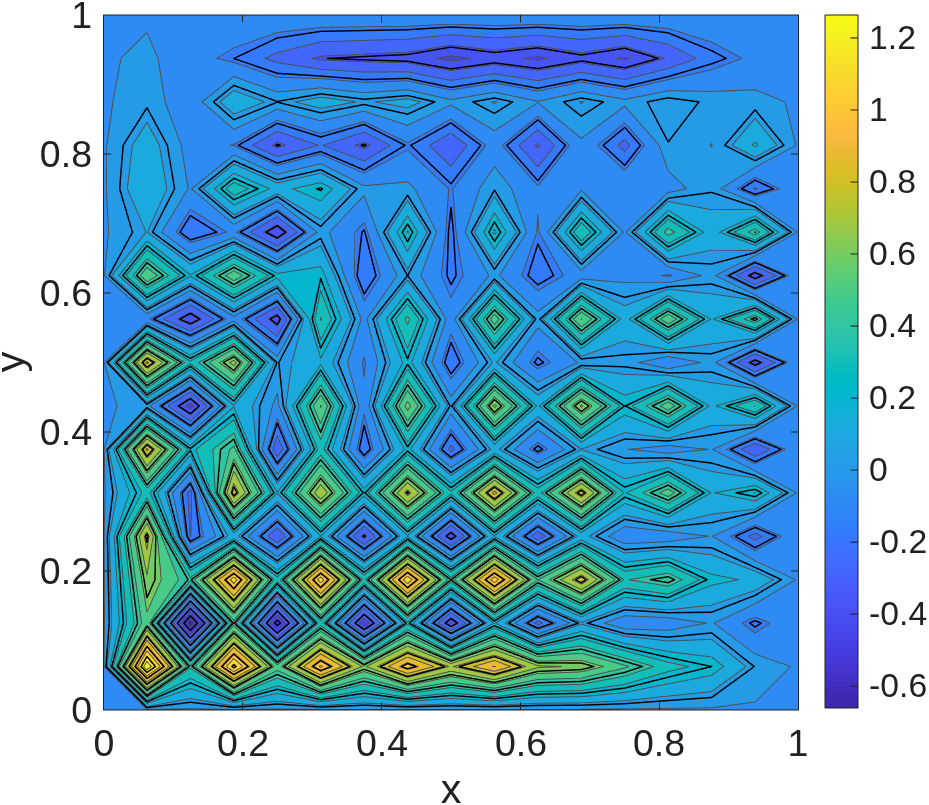}
        \caption{FEM solution at $z = 0.5$, $\alpha=10^{-3}$.}
        \label{fig:3d_bdy_layer_fem}
    \end{subfigure}
    \begin{subfigure}{0.45\linewidth}
        \centering
        \includegraphics[scale=0.28]{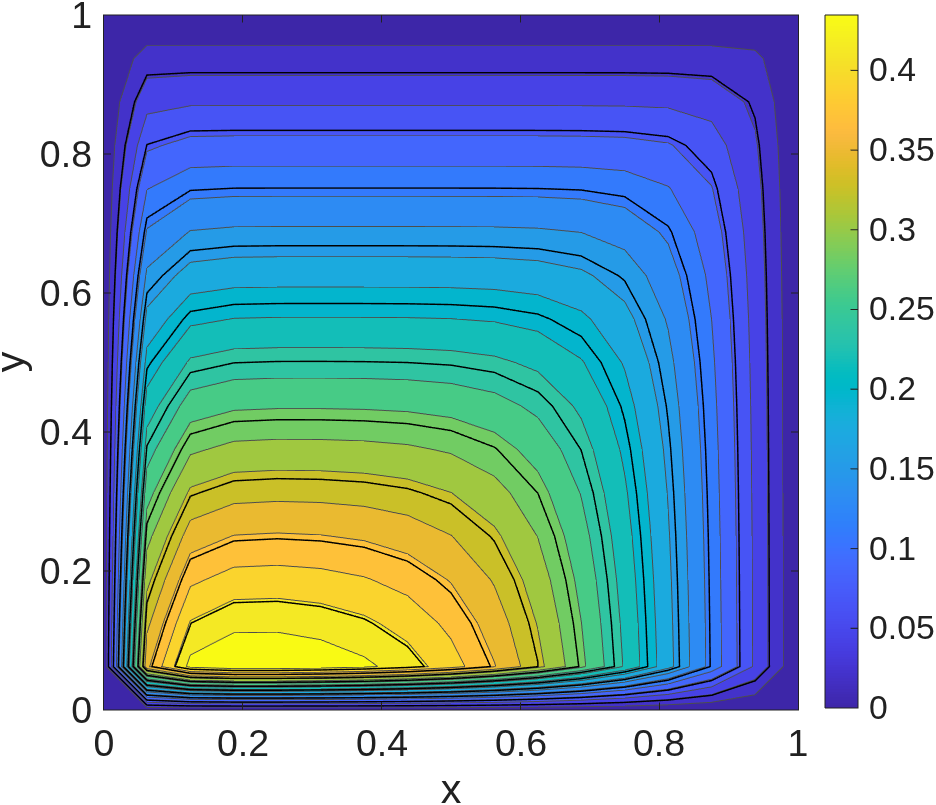}
        \caption{IBEAM solution at $z = 0.5$, $\alpha=10^{-3}$.}
        \label{fig:3d_bdy_layer_ibeam}
    \end{subfigure}
    \hskip 30pt
     \begin{subfigure}{0.45\linewidth}
        \centering
        \includegraphics[scale=0.28]{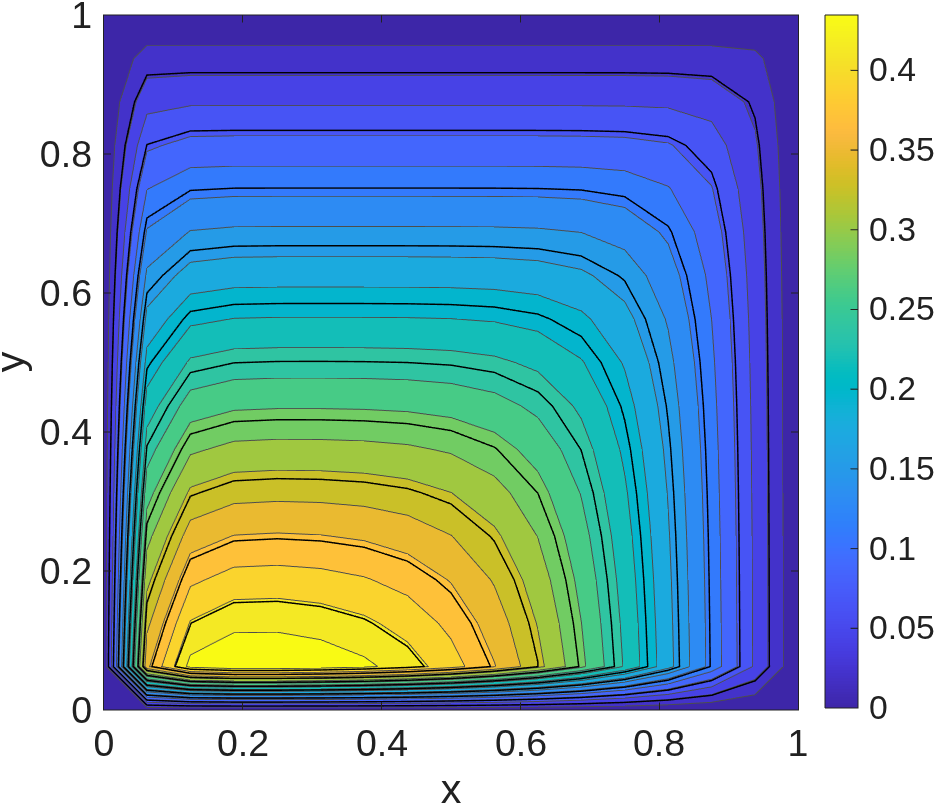}
        \caption{IBEAM solution at $z = 0.5$, $\alpha=10^{-6}$.}
        \label{fig:3d_bdy_layer_ibeam_6}
    \end{subfigure}
    \vskip -10pt
    \caption{Comparison of FEM and IBEAM solutions at the level curve $z = 0.5$ for $\alpha = 10^{-3}$, $\bm \beta = [1 \; 2 \; 0]^\top$,
    $f=1$, and $h=2^{-4}$, where IBEAM solution refers to the solution
    utilizing local stiffness matrix \eqref{ibeam_bi}. Unlike the FEM solution, which exhibits
    nonphysical oscillations, the IBEAM solution
    avoids oscillatory behavior, even for significantly smaller $\alpha=10^{-6}$.}
    \label{fig:3d_bdy_layer}
\end{figure}

\Cref{fig:3d_bdy_layer_ibeam,fig:3d_bdy_layer_ibeam_6} show a contour plot of the level curve at $z = 0.5$, where we observe a smooth solution unaffected by spurious oscillations; however, spurious oscillations appear in \Cref{fig:3d_bdy_layer_fem}, even for a larger value of $\alpha$.  Again, this validates the monotonicity properties established in \Cref{section:monotone}.

\subsubsection{Convergence Study}
Next, we perform a convergence study using the manufactured solution $u(x,y) = \sin(\pi x)\sin(\pi y)\sin(\pi z)$, with convection vector $\bm \beta = [0\;0\;1]^\top$, and the right-hand side chosen so that $u$ satisfies the convection-diffusion equation. We let the diffusion parameter take values $\alpha \in \{1, 0.1, 0.01, 0.001\}$ and consider meshes with mesh size $h = 2^{-2}$ to $h = 2^{-6}$.
For each mesh family discussed above, we again report the $L^2$ norm and $H^1$ semi-norm error at each refinement level.     We compare the $H^1$ semi-norm error results against $O(h)$ as established in \cite{eafe}.

\begin{figure}[h!]
    \centering
    \begin{subfigure}{0.45\linewidth}
        \centering
        \includegraphics[width=\linewidth]{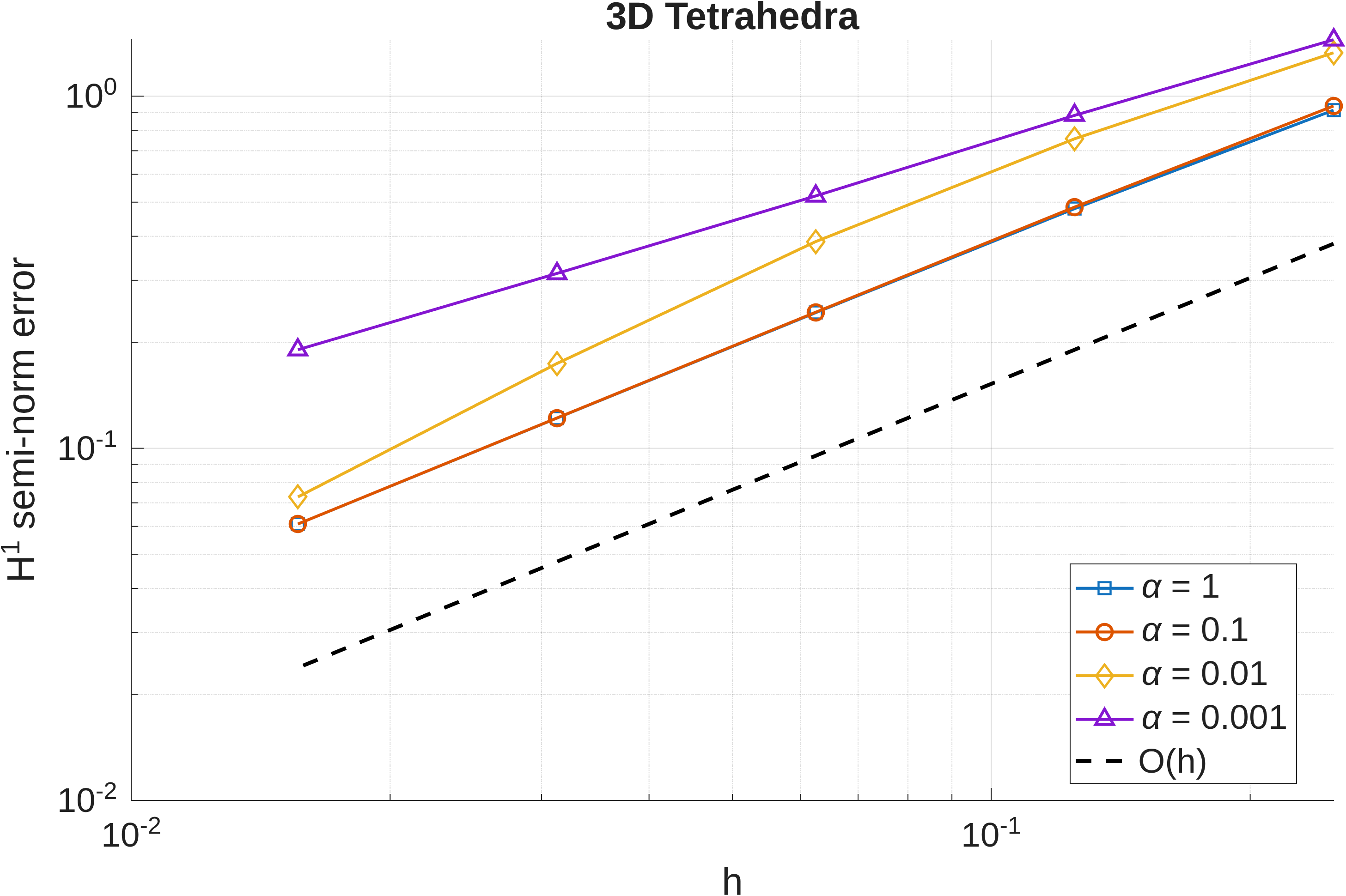}
        \label{fig:conv_3d_tet_h1}
    \end{subfigure}
    \hskip 30pt
    \begin{subfigure}{0.45\linewidth}
        \centering
        \includegraphics[width=\linewidth]{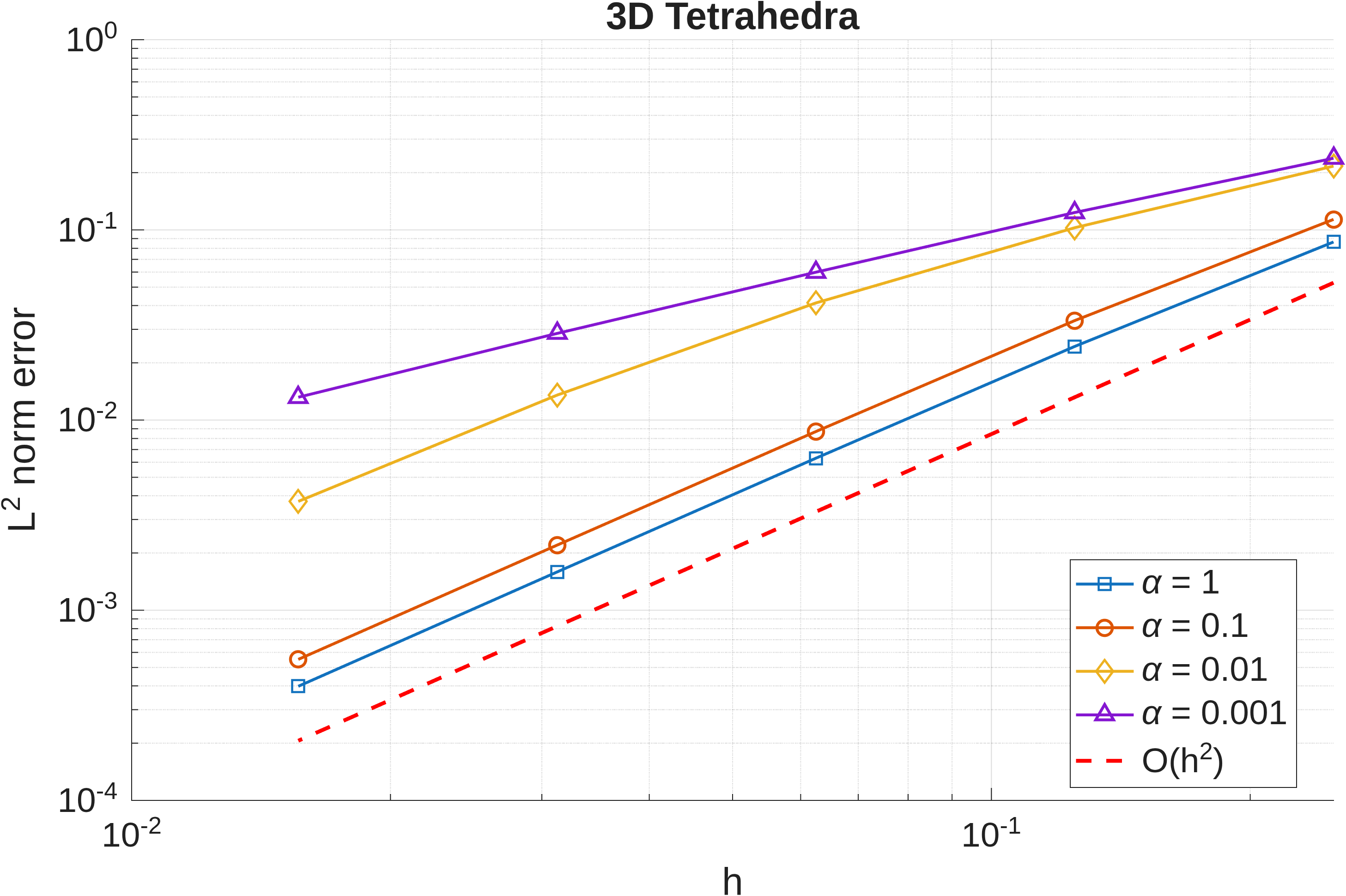}
        \label{fig:conv_3d_tet_l2}
    \end{subfigure}
    \vskip -10pt
    \caption{IBEAM error convergence rates on tetrahedral meshes for $\alpha \in \{1, 0.1, 0.01, 0.001\}$, plotted against mesh size $h$, for the $H^1$ semi-norm (left) and  $L^2$ norm (right) error.}
    \label{fig:conv_3d_tet}
\end{figure}

\begin{figure}[h!]
    \centering
    \begin{subfigure}{0.45\linewidth}
        \centering
        \includegraphics[width=\linewidth]{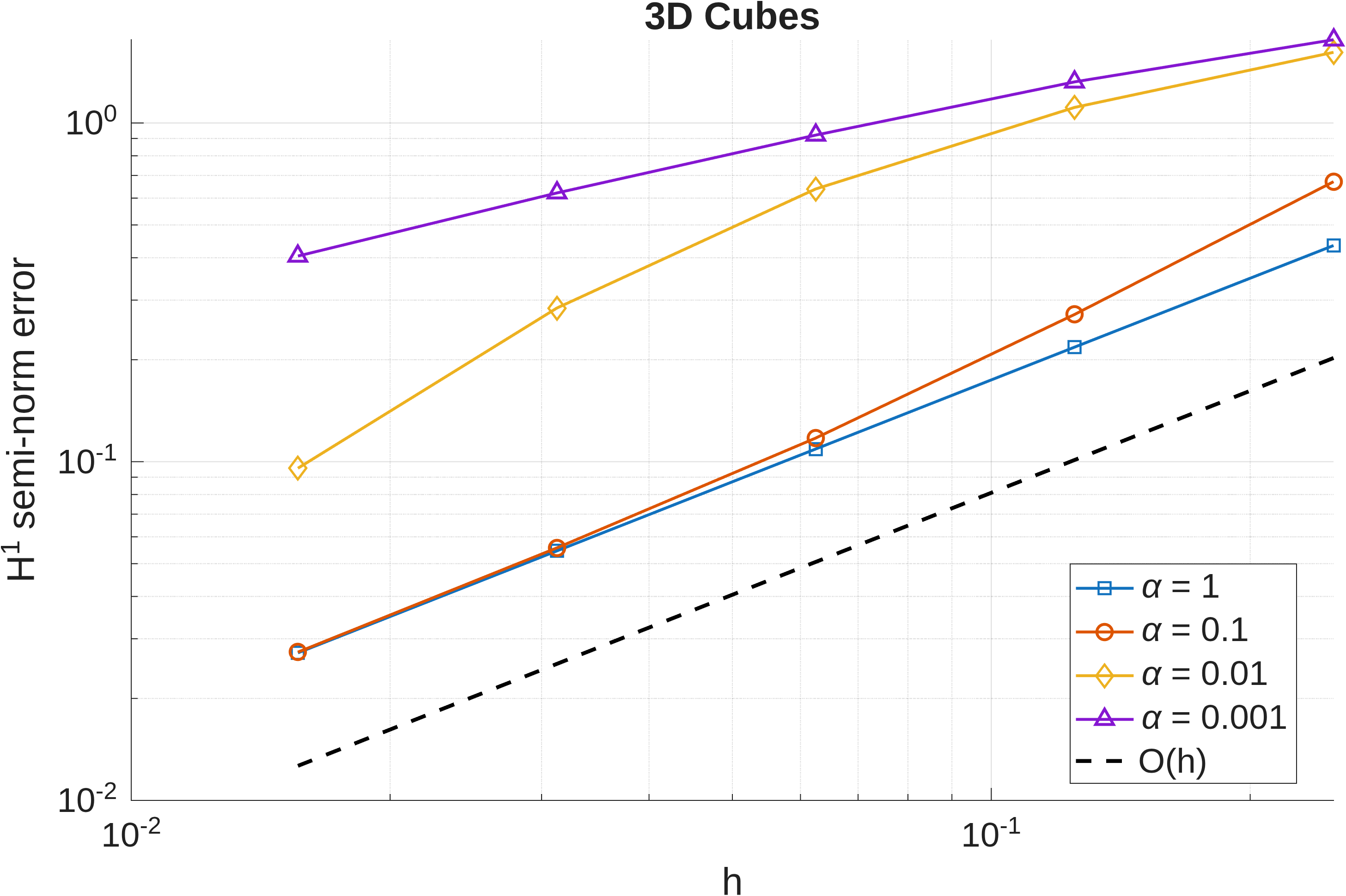}
        \label{fig:conv_3d_hex_h1}
    \end{subfigure}
    \hskip 30pt
    \begin{subfigure}{0.45\linewidth}
        \centering
        \includegraphics[width=\linewidth]{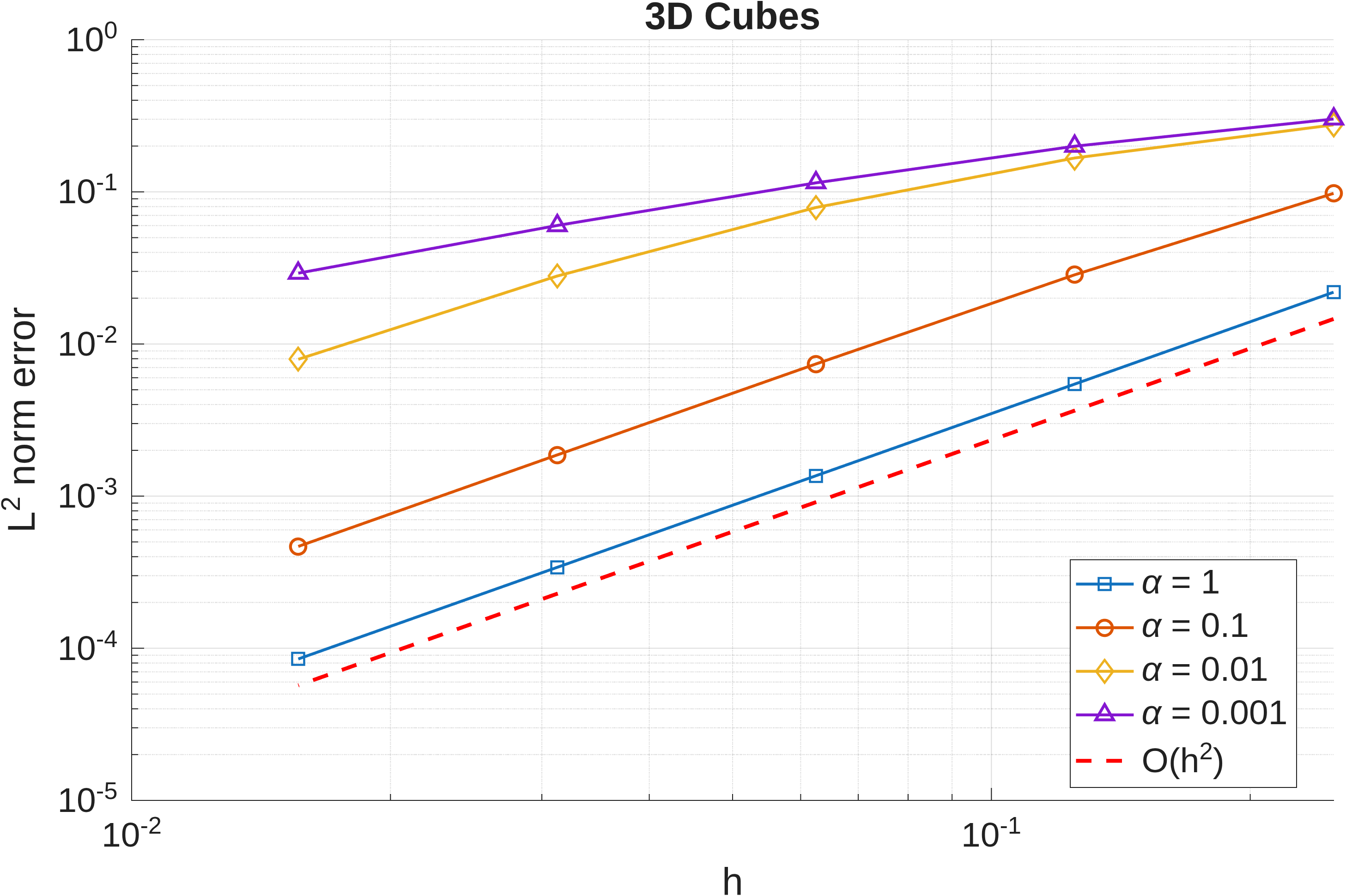}
        \label{fig:conv_3d_hex_l2}
    \end{subfigure}
    \vskip -10pt
    \caption{IBEAM error convergence rates on hexahedral meshes for $\alpha \in \{1, 0.1, 0.01, 0.001\}$, plotted against mesh size $h$, for the $H^1$ semi-norm (left) and  $L^2$ norm (right) error.}
    \label{fig:conv_3d_hex}
\end{figure}

On tetrahedral meshes (\Cref{fig:conv_3d_tet}), the computed rates agree with or exceed the established $O(h)$ rate in the $H^1$ semi-norm \cite{eafe} starting from the coarsest mesh for all but the smallest $\alpha$. We see that in the case of strongest convection dominance, the computed error rates in both the $H^1$ semi-norm and the $L^2$ norm fall short of the $O(h)$ and $O(h^2)$ bounds, respectively. We expect that further refinement, while maintaining monotonicity, ensures better convergence. Hexahedral meshes (\Cref{fig:conv_3d_hex}) demonstrate performance similar to that observed on tetrahedral meshes. These results show that IBEAM retains its predicted convergence behavior on non-simplicial three-dimensional meshes, provided sufficient refinement, which has not been addressed by EAFE or EAVE-type methods \cite{eave, safe, eafe}.

\section{Conclusion}\label{section:conclusion}

We developed an incidence-based edge-averaged method for scalar convection-diffusion problems on convex polytopal meshes in arbitrary dimensions. The proposed method extends the classical EAFE scheme through a graph-theoretic formulation based on vertex-edge incidence matrices. By introducing complementary edges and generalized barycentric coordinates, we established a matrix identity that relates the generalized Whitney edge mass matrix to a diagonal representation of the Poisson stiffness matrix, enabling an exponentially-fitted discretization on general convex polytopes. We also established sufficient conditions under which the resulting stiffness matrix is a nonsingular \(M\)-matrix, ensuring monotonicity and the discrete maximum principle. Numerical experiments on two- and three-dimensional meshes demonstrate first-order convergence in the \(H^1\) semi-norm and stable, nonoscillatory behavior in convection-dominated regimes. Approximately second-order convergence in the \(L^2\) norm is also observed in most cases.

One promising direction for future research is to apply the proposed method to time-dependent convection-diffusion problems through a space-time formulation \cite{adler2025unified}. By treating time as an additional coordinate, three-dimensional time-dependent problems can be formulated in four-dimensional space-time domains, providing an opportunity to investigate the numerical performance of IBEAM on four-dimensional polytopal meshes. Furthermore, the incidence-based framework suggests potential extensions to vector convection-diffusion problems in \(H(\mathrm{curl})\) and \(H(\mathrm{div})\), where edge-face and face-element incidence relations may be exploited to construct structure-preserving discretizations.

\section*{Acknowledgments}
The authors used AI to assist in generating portions of the code, in debugging, and scripting figure generation. All AI-assisted code was reviewed, tested, and validated by the authors, and the numerical results were verified independently. The authors assume full responsibility for all content. The authors also thank Ludmil Zikatanov for his valuable comments, suggestions, and insights toward the manuscript.

\bibliographystyle{siamplain}
\bibliography{references}
\end{document}